%% file: Arxiv_Submission0801.tex
\documentclass[times,sort&compress,3p]{elsarticleArx}
\journal{}
\usepackage[labelfont=bf]{caption}

\usepackage[T1]{fontenc}
\usepackage{amsmath,amsfonts,amssymb,amsthm,booktabs,color,epsfig,graphicx,hyperref,url,multirow,enumerate,mathtools,xcolor}

\theoremstyle{plain}
\newtheorem{thm}{Theorem}
\newtheorem{lma}{Lemma}

\newtheorem{cor}{Corollary}

\theoremstyle{definition}
\newtheorem{defn}{Definition}
\newtheorem{rem}{Remark}

\def\mc{\mathcal}
\def\T{\mc{T}}
\def\Tj{\mc{T}_j}
\def\Tk{\mc{T}_k}
\def\Rbb{\mathbb{R}}
\def\Nbb{\mathbb{N}}
\def\Cov{\mathrm{Cov}}
\def\Var{\mathrm{Var}}
\def\Pr{\mathrm{Pr}}
\def\tb{\mathbf{t}}
\def\inv{^{-1}}

\definecolor{DarkRed}{rgb}{.7,0,.4}

\newcommand{\normj}[1]{\left\lVert #1 \right\rVert_{j}} 
\newcommand{\normjk}[1]{\left\lVert #1 \right\rVert_{j, k}} 
\def\asympu{\overset{u}{\asymp}}
\def\asymplu{\overset{u}{\preceq}}

\def\Iijk{\mc{I}_{ijk}}

\newcommand{\pr}[1]{\left(#1\right)}

\def\muj{\mu_j}
\def\muk{\mu_k}
\def\hmuj{\hat{\mu}_j}
\def\hmuk{\hat{\mu}_k}
\def\hmujo{\hat{\mu}_{\textrm{obs},j}}
\def\hmujs{\hat{\mu}_{\textrm{subj},j}}
\def\gammajj{\gamma_{jj}}
\def\gammajk{\gamma_{jk}}
\def\hgammajk{\hat{\gamma}_{jk}}

\def\hgammajko{\hat{\gamma}_{\textrm{obs},jk}}
\def\hgammajks{\hat{\gamma}_{\text{subj},jk}}

\def\Uij{U_{ij}}
\def\Uik{U_{ik}}
\def\Uijl{U_{ij\ell}}
\def\Uijklm{U_{ijk\ell m}}
\def\Nij{N_{ij}}
\def\Nik{N_{ik}}

\def\ds{\mathrm{d}s}
\def\dt{\mathrm{d}t}
\def\du{\mathrm{d}u}
\def\dv{\mathrm{d}v}

\def\son{\sum_{i = 1}^n}

\def\soNij{\sum_{\ell = 1}^{N_{ij}}}

\def\sIijk{\sum_{(\ell,m) \in \Iijk}}
\def\i01{\int_0^1}
\def\iTj{\int_{\Tj}}
\def\iTjk{\int_{\Tj\times\Tk}}
\def\Njbar{\overline{N}_j}
\def\Nkbar{\overline{N}_k}
\def\Nbar{\overline{N}}
\def\NjH{N_j^H}
\def\NkH{N_k^H}
\def\NH{N^H}
\def\Njp{N_j^+}
\def\Np{N^+}

\def\Njtbar{\overline{N}_{j^2}}
\def\Ntbar{\overline{N}_{(2)}}
\def\Nmin{\overline{N}_{\textrm{min}}}
\def\Nmax{\overline{N}_{\textrm{max}}}
\def\NHmin{\NH_{\textrm{min}}}
\def\NHmax{\NH_{\textrm{max}}}
\def\Njkbar{\overline{N}_{jk}}
\def\Njkkbar{\overline{N}_{jk^2}}

\def\Njjkkbar{\overline{N}_{j^2k^2}}
\def\Njkp{N_{jk}^+}
\def\NjkH{N_{jk}^H}
\def\lsn{\limsup_{n \rightarrow \infty}}

\def\wjn{\overline{w}_{nj}}
\def\wjns{w_{nj}^*}
\def\vjkn{\overline{v}_{njk}}
\def\vjkns{v_{njk}^*}
\def\qjkn{q_{jkn}}
\def\Wij{W_{ij}}

\def\Wijk{W_{ijk}}

\def\bmuj{b_{\mu_j}}
\def\bgammaj{b_{\gamma_j}}
\def\bgammak{b_{\gamma_k}}
\def\bmu{b_\mu}
\def\bgamma{b_\gamma}
\newcommand{\Kb}[1]{K_{#1}}
\newcommand{\tKb}[1]{\tilde{K}_{#1}}
\def\bz{\beta_0}
\def\bo{\beta_1}
\def\bt{\beta_2}
\def\hbz{\hat{\beta}_0}
\def\hbo{\hat{\beta}_1}
\def\hbt{\hat{\beta}_2}

\def\Szj{S_{j0}}
\def\Soj{S_{j1}}
\def\Stj{S_{j2}}
\def\Rzj{R_{j0}}
\def\Roj{R_{j1}}

\def\Sjkzz{S_{jk00}}
\def\Sjkzo{S_{jk01}}
\def\Sjkoz{S_{jk10}}
\def\Sjkoo{S_{jk11}}
\def\Sjkzt{S_{jk02}}
\def\Sjktz{S_{jk20}}
\def\Sjkqr{S_{jkqr}}
\def\Rjkzz{R_{jk00}}
\def\Rjkzo{R_{jk01}}
\def\Rjkoz{R_{jk10}}
\def\pgjks{\frac{\partial \gamma_{jk}}{\partial s}}
\def\pgjkt{\frac{\partial \gamma_{jk}}{\partial t}}
\def\wij{w_{ij}}
\def\vijk{v_{ijk}}
\def\vijj{v_{ijj}}
\def\Nij{N_{ij}}
\def\Nik{N_{ik}}
\def\Yijl{Y_{ij\ell}}
\def\Yikm{Y_{ikm}}
\def\Zijklm{Z_{ijk\ell m}}
\def\Tijl{T_{ij\ell}}
\def\Tikm{T_{ikm}}
\def\Til{T_{i\ell}}

\def\epsijl{\epsilon_{ij\ell}}
\def\epsikm{\epsilon_{ikm}}

\DeclareFontEncoding{FMS}{}{}
\DeclareFontSubstitution{FMS}{futm}{m}{n}
\DeclareFontEncoding{FMX}{}{}
\DeclareFontSubstitution{FMX}{futm}{m}{n}
\DeclareSymbolFont{fouriersymbols}{FMS}{futm}{m}{n}
\DeclareSymbolFont{fourierlargesymbols}{FMX}{futm}{m}{n}
\DeclareMathDelimiter{\VERT}{\mathord}{fouriersymbols}{152}{fourierlargesymbols}{147}
\newcommand{\normsup}[1]{\left \VERT #1 \right \VERT} 
\newcommand{\normsupj}[1]{\left \VERT #1 \right \VERT_j} 
\newcommand{\normsupjk}[1]{\left \VERT #1 \right \VERT_{j,k}} 

\begin{document}

\begin{frontmatter}

\title{Mean and covariance estimation for discretely observed high-dimensional functional data: rates of convergence and division of observational regimes}

\author[1]{Alexander Petersen\corref{mycorrespondingauthor}}

\address[1]{Department of Statistics\\ 2152 WVB \\
Brigham Young University \\
Provo, UT 84602}

\cortext[mycorrespondingauthor]{Corresponding author. Email address: \url{petersen@stat.byu.edu}}

\begin{abstract}
Estimation of the mean and covariance parameters for functional data is a critical task, with local linear smoothing being a popular choice.  In recent years, many scientific domains are producing multivariate functional data for which $p$, the number of curves per subject, is often much larger than the sample size $n$.  In this setting of high-dimensional functional data, much of developed methodology relies on preliminary estimates of the unknown mean functions and the auto- and cross-covariance functions.  This paper investigates the convergence rates of local linear estimators in terms of the maximal error across components and pairs of components for mean and covariance functions, respectively, in both $L^2$ and uniform metrics.  The local linear estimators utilize a generic weighting scheme that can adjust for differing numbers of discrete observations $\Nij$ across curves $j$ and subjects $i$, where the $\Nij$ vary with $n$.  Particular attention is given to the equal weight per observation (OBS) and equal weight per subject (SUBJ) weighting schemes.  The theoretical results utilize novel applications of concentration inequalities for functional data and demonstrate that, similar to univariate functional data, the order of the $\Nij$ relative to $p$ and $n$ divides high-dimensional functional data into three regimes (sparse, dense, and ultra-dense), with the high-dimensional parametric convergence rate of $\left\{\log(p)/n\right\}^{1/2}$ being attainable in the latter two.  
\end{abstract}

\begin{keyword} 
Concentration inequalities \sep
High-dimensional data \sep
$L^2$ convergence \sep
Local linear smoothing \sep
Uniform Convergence
\MSC[2020] Primary 62G20 \sep
Secondary 62G08, 62G05
\end{keyword}

\end{frontmatter}

\section{Introduction\label{sec: intro}}

Over the last three decades, the foundational statistical principles underlying modern functional data analysis have been established \citep{rams:05,ferr:06,hsin:15,horv:12,koko:17}.  The relevant literature in more recent years suggests that the field is undergoing a transition, compelled by real data applications, in which the traditional univariate setting (one curve per subject) is being replaced with settings of larger and more complex structure, as evidenced by several review papers \citep{cuevas2014partial,goia2016introduction,mull:16:3,aneiros2019recent,aneiros2022functional}.
The authors of \cite{mull:16:3} described these new data sets as ``next-generation" functional data, including such cases as functions with non-linear constraints \citep{mull:16:1,hron:16,chen:21,ghod:22} or functions that take values in non-Euclidean spaces \cite{dube:20,dube:22}.  Similarly, \cite{kone:23} used the term ``second generation" to refer to large functional data sets with complex dependencies, including longitudinal \citep{grev:11,mull:12:3,park:15} and spatial \citep{deli:10,horm:11,zhan:22} functional observations, as well as multivariate functional data. Developments for multivariate functional data, in which observations across multiple curves are available for each subject, extend back to early methodological work on dimension reduction \citep{berr:11,mull:14:1,chio:14}, clustering \citep{jacq:14}, and regression models in which functions can appear as responses \citep{chio:16,zhu:17} or predictors \citep{wong:19}.  In the same way that classical multivariate data analysis led to modern developments in high-dimensional data analysis, some recent work in multivariate functional data methodology can best be described as high-dimensional functional data analysis, in which $p$, the number of curves for which observations are available per subject, is on the order of, or much larger than, the number of independent subjects $n$. Examples include interpretable dimension reduction \citep{hu:20} and discriminant analysis \citep{xue:22} for electroencephalography (EEG) data, and functional graphical models for both EEG \citep{zhu:16,qiao:19,sole:20,sole:22,lee:22,zhao:22} and functional magnetic resonance imaging (fMRI) data \citep{li:18,zapa:22}.  

A common element of the existing work on high-dimensional functional data is the requirement of preliminary estimates of the $p$ mean functions and the $p(p+1)/2$ distinct covariance functions.  The theoretical properties of these preliminary estimates determine those of the given method.  In terms of theoretical justification of their methodologies, previous works either assume that the curves are fully observed, or that they are observed on a regular, deterministic grid, with the number of observations per curve growing polynomially with the sample size.  The former scenario does not reflect the practical realities of functional datasets, while the latter does not cover important cases of designs in which observations can be random, irregular across the domain, heterogeneous across subjects, or any combination of these; moreover,  the number of observations per curve per subject may not diverge with $n$, or may do so at application-specific rates. For instance, \cite{zhon:19,fang:20} illustrate methodologies using biomedical data that are described therein as ``high-dimensional longitudinal data."  As longitudinal data can be modeled as functional data with relatively few temporal observations compared to the number of independent sampling units, such data constitute an example of sparsely observed high-dimensional functional data for which crucial theory has not yet been developed.  It should be emphasized that, while a multivariate functional data set with a small number of curves $p$, but a large number of observations per curve relative to the sample size $n$, may accurately be described as high-dimensional, this is not the sense in which this term is used in this paper. Rather, dimensionality refers to the number of curves per subject, and will always be assumed to be high, while the number of discrete observations per curve ranges from sparse to dense, and determines the convergence rate of the estimators.

The contribution of this paper is to derive rates of convergence for nonparametric estimators of the mean and covariance functions for high-dimensional functional data.  Specifically, local linear estimators will be studied in terms of their accuracy in both the $L^2$ (in probability) and uniform (almost surely) metrics.  This choice was made due to critical previous work in the setting of univariate or ``first-generation" functional data \citep{li:10,zhan:16}, corresponding to a special case of the developments herein when $p = 1.$  Both of these works studied local linear estimators of mean and covariance functions in the manner described, the latter in a more general fashion.  In particular, \cite{li:10} applied a specific weighting scheme to the observations in the estimation criterion, whereas \cite{zhan:16} considered a general weighting scheme and also demonstrated pointwise asymptotic normality of the estimators.  Additionally, as the covariance estimates are slightly different in these two papers, that of \cite{zhan:16} will be used in this paper; for a detailed and intuitive justification of this choice, the reader is referred to Section 6 of \cite{zhan:16}.

As the number of observations from a given curve can vary across subjects and across different components of the multivariate functional data, and since these observations are inherently correlated, it is reasonable to assign different weights for observations coming from distinct subjects or observational units in the data set.  Two common choices are the subject (SUBJ) scheme, in which the total weight of all observations from a given subject is the same across subjects, and the observation (OBS) scheme, in which all observations share the same weight. 
The former choice was used in \cite{li:10}, while \cite{zhan:16} studied a generic weighting scheme with specific focus on the SUBJ and OBS schemes.  The following two critical insights were provided by \cite{zhan:16}.  First, if the number of observations for a given component curve is highly heterogeneous across subjects, the OBS scheme can perform poorly due to higher weight being granted to a group of potentially highly correlated observations; on the other hand, if they are homogeneous in an explicit sense, the OBS scheme will yield a rate no worse than the SUBJ scheme.  Second, depending on the behavior of the average (for the OBS scheme) or hyperbolic average (for the SUBJ scheme) number of observations per subject relative to $n$, prescriptive bandwidth choices divide functional data sets into three observational regimes: sparse (or non-dense), dense, and ultra-dense.  In the sparse scenario, the rate can be anywhere between the nonparametric (inclusive) rates of $n^{-2/5}$ and $n^{-1/3}$ for mean and covariance estimation, respectively, and the parametric rate of $n^{-1/2}$ (exclusive).  For dense and ultra-dense data, the parametric rate is always attainable; the distinguishing feature of ultra-dense data is that the bias decays at a faster rate than the stochastic error, whereas these two are matched for dense data.

The results derived in this paper generalize the above findings to the high-dimensional setting where $p$ diverges with $n$.  Specifically, a comparable rate of convergence for the largest error across different mean or covariance estimates is derived in the high-dimensional setting, with the only differences being an inflation of the stochastic rate by $\{\log(p)\}^{1/2}$ and an additional term arising from the involvement of higher-order moments in the newly derived bounds; see Theorems~\ref{thm: meanGen}--\ref{thm: CovSupGen}.  The OBS scheme is again found to be adversely affected by heterogeneous numbers of observations, where the quantification of homogeneity (see \eqref{eq: OBShom} and \eqref{eq: OBShom2}) is again explicit but more strict than that of \cite{zhan:16}.  On the other hand, if the homogeneity condition is satisfied, the OBS scheme is never worse than the SUBJ scheme; see Corollaries~\ref{cor: meanObsSubjL2}, \ref{cor: meanObsSubjSup}, \ref{cor: CovObsSubjL2}, and \ref{cor: CovObsSubjSup}. For these two schemes, optimal bandwidth decay rates are also prescribed, yielding observational regime divisions analogous to those of \cite{zhan:16} in the case of mean estimation.  For covariance estimation, while the parametric rate $\left\{\log(p)/n\right\}^{1/2}$ is still attainable for dense and ultra-dense data, the results in this paper do not distinguish between these in terms of bias decay; see Remark~\ref{rem: CovUD} for a more detailed explanation.

The remainder of the paper is organized as follows.  Section~\ref{sec: methods} provides definitions of the functional targets and estimators, as well as two particular classes of observational designs that will be considered.  Sections~\ref{sec: mean} and \ref{sec: cov} provide all technical assumptions and theoretical results corresponding to the mean and covariance estimators, respectively.  The paper concludes with a brief discussion in Section~\ref{sec: disc}, while all proofs are provided in Section~\ref{sec: App}.  The supplementary material includes simulations that illustrate the performance of the OBS and SUBJ schemes for different observational designs in the high-dimensional regime $p > n,$ including a discussion of computational aspects and challenges of constructing a large number of smoothing estimators.

\section{Methodology\label{sec: methods}}

For $p \in \Nbb,$ let $\Tj$, $j \in \{1,\ldots,p\},$ be compact intervals of the real line and $\mc{T}^p = \bigtimes_{j = 1}^p \Tj$ their Cartesian product.  It will be assumed throughout that $p$ diverges with $n$ such that $\log(p)/n \rightarrow 0.$  Let $\{X(\tb) \in \Rbb^p;\, \tb = (t_1,\ldots,t_p) \in \mc{T}^p\}$ be a multivariate $L^2$ stochastic process, that is, $X(\tb) = (X_1(t_1),\ldots,X_p(t_p))^\intercal$ satisfies $E(X_j^2(t)) < \infty$ for all $t \in \Tj$ and all $j \in \{1,\ldots,p\}.$  The primary population targets for which estimates are typically sought in functional data analysis are the mean and covariance functions of $X$, which will be denoted in this paper by
\begin{equation}
\label{eq: MeanCov}
\muj(s) = E[X_j(s)], \quad \gammajk(s, t) = \Cov[X_j(s), X_k(t)],
\end{equation}
for $(s, t) \in \Tj \times \Tk$, $j,k \in \{1,\ldots,p\}.$  When $j = k,$ $\gamma_{jj}$ is referred to as the $j$-th auto-covariance function of $X$, reflecting the intracurve dependence for a given component function, while $\gammajk$ for $j \neq k$ are the cross-covariance functions corresponding to dependence between two distinct functions or curves.  Note that generic arguments $s$ and $t$ are used for all of these functions regardless of the domain; these arguments will be referred to as time points, although the functional domains $\Tj$ need not correspond to time.  The ultimate aim of this paper is to determine rates of convergence for estimators of these targets constructed from a suitable sample in the high-dimensional regime.

Suppose $X_1,\ldots,X_n$ are independently and identically distributed as $X,$ and write each of these independent multivariate processes as $X_i(\tb) = (X_{i1}(t_1),\ldots,X_{ip}(t_p))^\intercal, $ $i \in \{1, \ldots, n\}.$  A key challenge for functional data is that the processes $X_i$ are never fully observed along the continua; rather, the collected data correspond to a finite number of measurements over a grid of points for each curve.  The scheme by which the grid points arise and how these are modeled are referred to as the observational design.  Two design settings, described in Section~\ref{ss: designs}, will be considered in this paper as extensions of the random design setting treated in previous papers on the topic for univariate functional data \citep{li:10, zhan:16} in which $p = 1$.  In both settings, data for each curve are collected at a random collection of points along the domain, where the number of points per curve may be bounded or diverge with the sample size $n$.  

Prior to considering the particular nuances of the two designs, consider a general model and corresponding mean and covariance estimators. For each subject $i$ and component $j$, the observed data are modeled as
\begin{equation}
\label{eq: DataModel}
\Yijl = X_{ij}(\Tijl) + \epsijl = \mu_j(\Tijl) + \Uij(\Tijl) + \epsijl, \quad \ell \in \{1,\ldots,N_{ij}\},
\end{equation}
where $\Uij = X_{ij} - \muj$ is the centered process, $N_{ij}$ are the number of observation points $\Tijl$ for this curve, and $\epsijl$ are error variables with mean zero, independent of the processes $\Uij;$ interdependence of the errors will be left unspecified for the moment as it depends on the observational design. Throughout, the $N_{ij}$ will be considered as deterministic, though varying with $n$, and their behavior as $n$ diverges will heavily impact the convergence rates. 

Following previous work \citep{li:10,zhan:16},  local linear techniques will be used to estimate the functions in \eqref{eq: MeanCov}.  Let $K$ be a univariate probability density function and, for a bandwidth $b > 0, $ define $\Kb{b}(\cdot) = b\inv K(\cdot/b).$  For each $j \in \{1,\ldots,p\}$, let $\bmuj > 0$ and consider positive weights $w_{ij}$ that satisfy $\son w_{ij}N_{ij} = 1.$  For $t \in \Tj,$ define $\hmuj(t) = \hbz,$ where
\begin{equation}
\label{eq: MeanEst}
(\hbz, \hbo) = \arg \min_{\bz, \bo} \son w_{ij} \soNij \Kb{\bmuj}(\Tijl - t)\left\{\Yijl - \bz - \bo(\Tijl - t)\right\}^2.
\end{equation}

Covariance estimation then follows by defining the raw covariance terms $\Zijklm = \{\Yijl - \hmuj(\Tijl)\}\{Y_{ikm} - \hmuk(\Tikm)\}$ for any $j,k \in \{1,\ldots,p\}$, $\ell \in \{1,\ldots N_{ij}\}$, and $m \in \{1,\ldots,N_{ik}\}.$  To avoid inducing bias in the covariance estimation, raw covariances are typically removed from the estimation if they contain dependencies between the noise variables.  To maintain generality for the moment, let $\Iijk \subset \{1,\ldots,\Nij\} \times \{1,\ldots,\Nik\}$ denote a suitable subset of index pairs $(\ell,m)$ for which the raw covariances will be included in the estimation; further specification of these will be given in Section~\ref{ss: designs}.  Let $\bgammaj,\bgammak > 0$ and consider positive weights $\vijk$ satisfying $\son \vijk|\Iijk| =  1.$  For $(s, t) \in \Tj \times \Tk,$ define
\begin{equation}
\label{eq: CovEst}
(\hbz,\hbo,\hbt) = \arg \min_{\bz, \bo, \bt} \son v_{ijk} \sIijk \Kb{\bgammaj}(\Tijl - s)\Kb{\bgammak}(\Tikm - t)\left\{\Zijklm - \bz - \bo(\Tijl - s) - \bt(\Tikm - t)\right\}^2.
\end{equation}
Then the estimate of $\gammajk(s,t)$ is $\hgammajk(s, t) = \hbz.$  For simplicity, the same bandwidth $\bgammaj$ is used for estimating the auto-covariance $\gammajj$ and any cross-covariance $\gammajk.$  For the kernel $K$, the following assumptions are required.
\begin{itemize}
    \item[A1] $K$ is a probability density function with support $[-1,1],$ is symmetric about zero, and is of bounded variation.
    \item[A2] $K$ is Lipschitz continuous.
\end{itemize}
Assumption A1 is ubiquitous in the kernel smoothing literature, while assumption A2 allows for simplification of the proofs of results involving the uniform metric as it governs the smoothness of the estimators \citep{zhan:16}.  While not strictly necessary, assumption A2 is not restrictive in practice, as it is satisfied for commonly used kernels, such as the Epanichnikov and Gaussian kernels.  However, it can be omitted at the cost of more cumbersome arguments \citep{li:10}.

The estimators defined by \eqref{eq: MeanEst} and \eqref{eq: CovEst} will be assessed in terms of their convergence rates in the $L^2$ and uniform norms, denoted by $\normj{f} = \left\{ \iTj f^2(t) \dt\right\}^{1/2}$, $\normsupj{f} = \sup_{t \in \Tj} |f(t)|$, $\normjk{g} = \left\{ \iTjk g^2(s,t) \ds \dt\right\}^{1/2}$, and $\normsupjk{g} = \sup_{(s,t) \in \Tj \times \Tk} |g(s,t)|.$
Specifically, rates of the convergence for the maximal of these norms across $j$ or $(j, k)$ will be determined that will allow for their consistent estimation so long as $\log(p)$ grows more slowly than the sample size.  As in \cite{zhan:16}, rates will be determined for generic weighting schemes, and special attention will be given to the so-called observation (OBS) and subject (SUBJ) weighting schemes.  Define $\Njbar = n\inv\son N_{ij}$ as the average number of observations of the $j$-th curve across subjects.  Then the OBS weights are $\wij = 1/(n\Njbar)$ and $\vijk = 1/(\son |\Iijk|)$.  For the SUBJ scheme, one has $\wij = 1/(n\Nij)$ and $\vijk = 1/(n|\Iijk|)$.  Write $\hmujo$ and $\hgammajko$ for the mean and covariance estimators under the OBS scheme, and $\hmujs$ and $\hgammajks$ for those under the SUBJ scheme.

\subsection{Observational Designs\label{ss: designs}}

The first design setting treats the general case in which the observational time points may be different in number and location across curves.  This will be referred to as the fully random design.
\begin{defn}
\label{def: FRdesign}
The observation times follow a fully random (FR) design if, for each $j \in \{1,\ldots,p\},$ $\Tijl, $ $\ell \in \{1,\ldots, \Nij\}$ and $i \in \{1, \ldots, n\},$ are independently distributed on $\Tj$ with probability density $f_j$ and are independent across $j$.
\end{defn}
When the data come from an FR design, the errors $\epsijl$ will be assumed independent across $i$ and $\ell$ and, for covariance estimation, they will also be assumed independent across $j$.  When estimating covariance functions, raw covariances are thus only excluded from \eqref{eq: CovEst} when $j = k$ and $\ell = m$.  Thus, when an FR design is assumed, one has
\begin{equation*}
\Iijk = \left\{ \begin{array}{ll} \left\{(\ell,m):\, \ell \in \{1,\ldots,\Nij\},\, m \in \{1,\ldots,N_{ik}\}\right\}, & j \neq k, \\
\left\{(\ell, m):\, \ell,m \in \{1,\ldots,\Nij\} \textrm{ and } \ell \neq m\right\}, & j = k.\end{array} \right.
\end{equation*}

While some multivariate functional data sets are most appropriately modelled with an FR design, many have a specialized structure, especially when considering examples where the multivariate dimension $p$ is large.  In these cases, all functions share the same domain and, for each subject, all $p$ curves are observed simultaneously along a common set of timepoints.  This will be referred to as the simultaneous random design. 
\begin{defn}
\label{def: SRdesign}
The observation times follow a simultaneous random (SR) design if, for each $j \in \{1,\ldots,p\},$ $\Tj = \T$ for a compact interval $\T$, $N_{ij} = N_i,$ and $\Tijl = \Til,$ where $\Til$,  $\ell \in \{1,\ldots, N_i\},$ and $i \in \{1, \ldots, n\},$ are independently distributed on $\T$ with probability density $f.$
\end{defn}
When the data are assumed to come from an SR design, error variables $\epsijl$ and $\epsilon_{ik\ell}$ that are measured at a common timepoint $\Til$ may be dependent. Thus, under this design, the index sets for raw covariances will take the form
\begin{equation*}
    \Iijk = \mc{I}_i = \left\{(\ell, m):\, \ell,m \in \{1,\ldots,N_i\} \textrm{ and } \ell \neq m\right\}.
\end{equation*}
The following notational simplifications will also be assumed under the SR design.  The weights $\wij$ and $\vijk$ will be assumed independent of $j$ and $k$, denoted as $w_i$ and $v_i$ when appropriate.  In the OBS scheme, these take the form $w_i = 1/(n\Nbar)$ and $v_i = 1/(\son N_i(N_i - 1)),$ with $\Nbar = n\inv \son N_i;$  in the SUBJ scheme, they are $w_i = 1/(nN_i)$ and $v_i = 1/(nN_i(N_i - 1)).$  In addition, it will be assumed that a single bandwidth $\bmu$ is used to estimate all $p$ mean functions, and a common bandwidth $\bgamma$ is used for all $p(p+1)/2$ covariance functions.

\subsection{Technical Assumptions on Model Parameters}

The following assumptions on model \eqref{eq: DataModel} apply to all relevant theoretical results in Sections~\ref{sec: mean} and \ref{sec: cov}.
\begin{itemize}

\item[B1] The collections of functional data $\{X_i:\, i \in \{1, \ldots, n\}\}$, observation times $\{\Tijl:\, \ell \in \{1,\ldots, \Nij\},\, j \in \{1,\ldots,p\},\, i \in \{1, \ldots, n\}\}$, and errors $\{\epsijl:\, \ell \in \{1,\ldots, \Nij\},\, j \in \{1,\ldots,p\},\, i \in \{1, \ldots, n\}\}$ are independent of each other.  Furthermore, the $X_i$ are iid across $i$, and the data arise from either the FR design in Definition~\ref{def: FRdesign} or the SR design in Definition~\ref{def: SRdesign}.  Lastly, with $|\Tj|$ denoting the length of the interval $\Tj,$ $0 < \lim_{n \rightarrow \infty} \min_{j \in \{1,\ldots,p\}} |\Tj| \leq \lim_{n \rightarrow\infty} \max_{j \in \{1,\ldots,p\}} |\Tj| < \infty.$

\item[B2] Under an FR design, for each $j \in \{1,\ldots,p\},$ $f_j$ is a twice differentiable probability density function on $\Tj.$  Moreover, defining $m = \lim_{n\rightarrow \infty} \min_{j \in \{1,\ldots,p\}}\inf_{t \in\Tj} f_j(t)$ and $M = \lim_{n \rightarrow \infty} \max_{j \in \{1,\ldots,p\}} \max\left\{\normsupj{f_j},\normsupj{f_j''}\right\},$ $0 < m < M < \infty.$
Under an SR design, the preceding holds with $f_j \equiv f.$

\item[B3] The $\mu_j$ are twice differentiable, and $\lim_{n \rightarrow \infty} \max_{j \in \{1,\ldots,p\}} \normsupj{\mu_j''} < \infty.$

\item[B4] The $\gammajk$ are twice partially differentiable, and 
$$
\lim_{n \rightarrow \infty} \max_{j,k \in \{1,\ldots,p\}} \max\left(\normsupjk{\frac{\partial^2 \gammajk}{\partial s^2}}, \normsupjk{\frac{\partial^2 \gammajk}{\partial s \partial t}}, \normsupjk{\frac{\partial^2 \gammajk}{\partial t^2}}\right) < \infty.
$$
\end{itemize}

Assumption B1 stipulates the independence of the different random components in the model and the nature of the observation times; to simplify later conditions on bandwidths, it also asserts that the size of the domains for the different functional data components are of the same order.  Assumptions B2--B4 are regularly assumed in the case of univariate functional data, and are strengthened here in order for bounds to be uniformly controlled as $p$ diverges.

\section{Mean Estimation\label{sec: mean}}

In this section, rates of convergence will be provided for mean estimation in the high-dimensional regime for a generic weighting scheme under both FR and SR designs.  These results demonstrate that the effects of high-dimensionality are as expected, in that the effective sample size is reduced to $n/\log(p).$  Hence, in determining the divisions between non-dense, dense, and ultra-dense functional data \citep{zhan:16}, the rate of $\{\log(p)/n\}^{1/2}$ will be critical.

\subsection{\texorpdfstring{$L^2$}{L2} Convergence\label{ss: L2mean}}

Consider the estimator defined in \eqref{eq: MeanEst} for general sequences of weights $\wij.$  Let $\overline{w}_{nj} = \max_{i \in \{1,\ldots,n\}} w_{ij}N_{ij}$. The following assumptions on the distributional characteristics of the functional observations and the asymptotic behavior of bandwidths are required.  Let $p \vee n$ denote the maximum of $p$ and $n$. 
\begin{itemize}

    \item[C1] For each $j \in \{1,\ldots,p\}$ and $t \in \Tj,$ $\Uij(t)$ is sub-Exponential with parameter $\theta_j(t) > 0$, that is, $E[\exp\{\lambda \Uij(t)\}] \leq \exp\{\lambda^2\theta_j^2(t)/2\}$ for all $|\lambda| < \{\theta_j(t)\}\inv.$  Furthermore, 
    $
    \overline{\theta} = \lim_{n \rightarrow \infty} \max_{j \in \{1,\ldots,p\}} \sup_{t \in \Tj} \theta_j(t) < \infty.
    $
    
    \item[C2] For each $j$, the errors $\epsijl$ are sub-Exponential random variables with parameter $\sigma_j,$ that is, $E[\exp\{\lambda \epsilon_{ijl}\}] \leq \exp\{\lambda^2\sigma_j^2/2\}$ for all $|\lambda| < \sigma_j\inv.$  Furthermore,
    $
    \sigma = \lim_{n\rightarrow \infty} \max_{j \in \{1,\ldots,p\}} \sigma_j < \infty.
    $
    In addition, under either an FR or SR design, the $\epsijl$ are independent across $i$ and $\ell$ for any fixed $j$, but may be dependent across $j$.
    
    \item[C3] The bandwidths satisfy $\log(p)\max_{j \in \{1,\ldots,p\}} \son \wij^2\Nij(\bmuj\inv + \Nij - 1)  \rightarrow 0$, \mbox{$\log(p)\max_{j \in \{1,\ldots,p\}} \bmuj\inv\wjn \rightarrow 0$}, and $\max_{j \in \{1,\ldots,p\}} \bmuj \rightarrow 0$. In addition, $\log(p)\max_{j \in \{1,\ldots,p\}} \bmuj^{-2} \son \wij^2 \Nij \rightarrow 0$ and $\bmu^{-2}\son w_i^2N_i \rightarrow 0$ under the FR and SR designs, respectively.
    
    \item[C4] There exists $\alpha > 0$ such that 
    $$
    \frac{n^{\alpha}\max_{j \in \{1,\ldots,p\}} \left[\left\{\log(p\vee n)\son \wij^2\Nij(\bmuj\inv + \Nij - 1)\right\}^{1/2} + \log(p\vee n)\bmuj\inv\wjn\right]}{\max_{j \in \{1,\ldots,p\}} \bmuj^{-2}} \rightarrow \infty.
    $$
\end{itemize}

Assumptions C1 and C2 are stronger tail conditions compared to previous work for the case $p = 1$, in which only moment bounds are used to derive concentration inequalities. In ordinary (non-functional) high-dimensional data analysis, two common classes of tail behavior for analyzing mean estimation are polynomial and sub-Exponential tails.  The latter lead to probability bounds with exponential decay, and thus faster rates of convergence when $p$ diverges.  While acknowledging that rates and division of observational regimes for different classes of tail behavior are certainly of interest, this paper focuses on sub-Exponential tails as in assumptions C1 and C2, leaving other cases for future work. 
Assumptions C1 and C2, in conjunction with the stipulations on the bandwidths in assumption C3, thus strengthen the second moment assumptions made in \cite{zhan:16} in order to obtain $L^2$ convergence results in Section~\ref{ss: L2mean} below that are comparable to previous work. In the case of uniform convergence, assumption C1 is not directly comparable to its counterpart in \cite{li:10,zhan:16}, which require that $\normsupj{X_{ij} - \mu_j}$ have bounded $r$-th moment for some $r > 2.$  This latter condition can be viewed as a smoothness assumption on the sample paths of the $X_{ij}$, and does not imply, neither is it implied by, assumption C1 above.  Nevertheless, as will be seen in Section~\ref{ss: Supmean}, the same rates of convergence are also obtained for the uniform metric under tail conditions of assumptions C1 and C2 and under the conditions specified by assumptions C3 and C4 on the bandwidths.  The proofs of all results in the section can be found in Section~\ref{AppC}, with auxiliary lemmas and their proofs given in Section~\ref{AppB}.

\begin{thm}
\label{thm: meanGen}
Under assumptions A1, B1--B3, and C1--C3, 
\[
\max_{j \in \{1,\ldots,p\}} \normj{\hmuj - \muj} = O_P\left(\max_{j \in \{1,\ldots,p\}} \left[\bmuj^2 + \left\{\log(p)\son\wij^2\Nij(\bgammaj\inv + \Nij - 1)\right\}^{1/2} + \log(p)\bmuj\inv\wjn \right]\right).
\]
\end{thm}

\begin{rem}
    \label{rem: design}
   The rate given in Theorem~\ref{thm: meanGen} applies to both FR and SR designs. 
  Under the latter design, the maximum over $j$ in the rate is redundant, requiring less stringent requirements on the bandwidth; see Corollary~\ref{cor: meanOptL2_sim}.
\end{rem}

\begin{rem}
\label{rem: rateComp}
The rate in Theorem~\ref{thm: meanGen} is obtained by establishing an exponential tail bound on deviation probabilities for each mean estimate, followed by application of the union bound, which leads to the appearance of $\log(p)$ in last two terms making up the stochastic part of the rate; the bias is not affected.  In order to compare this result with the rate given in \cite{zhan:16}, the arguments in the proof of Theorem~\ref{thm: meanGen} can be applied to a single mean estimate to yield
$$
\normj{\hmuj - \muj} = O_P\left[\bmuj^2 + \left\{\son\wij^2\Nij(\bgammaj\inv + \Nij - 1)\right\}^{1/2} + \bmuj\inv\wjn\right],
$$
while the rate of \cite{zhan:16} omits the final term. Its appearance in the rate derived in this paper can be explained as follows. The result relies on an exponential tail bound obtained by applying Theorem 2.5 of \cite{bosq:00} and involves all pointwise moments of the functional data as well as the error moments, as opposed to only the first two moments used by \cite{zhan:16} that were sufficient to obtain the $L^2$ convergence rate in the setting $p = 1$ using Chebyshev's inequality.  Specifically, to leverage Theorem 2.5 of \cite{bosq:00}, for any $\nu \geq 2$, one must obtain a moment bound of the form
$$
\son \wij^\nu E\left[\left\{\soNij \Kb{\bmuj}(\Tijl - t)\right\}^\nu\right] \leq c_{n1}^2c_{n2}^{\nu - 2}.
$$
As $\soNij \Kb{\bmuj}(\Tijl -t)$ is almost surely bounded above by a multiple of $\bmuj\inv\Nij$ and $E\left[\left\{\soNij \Kb{\bmuj}(\Tijl - t)\right\}^2\right]$ is bounded by a multiple of $\Nij(\bmuj\inv + \Nij - 1)$, for some $C > 0,$ 
\[
\begin{split}
\son \wij^\nu E\left[\left\{\soNij \Kb{\bmuj}(\Tijl - t)\right\}^\nu\right] &\leq 
(C\bmuj\inv)^{\nu - 2} \left\{C\son \wij^\nu \Nij^{\nu - 1}(\bmuj\inv + \Nij - 1)\right\} \\
&\leq \left\{C \son \wij^2\Nij(\bmuj\inv + \Nij - 1)\right\}(C\bmuj\inv\wjn)^{\nu - 2}
\end{split}
\]
Hence, the final term in the rate cannot be eliminated under general bandwidth sequences. Nevertheless, when applied to the OBS and SUBJ weighting schemes, Corollaries~\ref{cor: meanOptL2} and \ref{cor: meanOptL2_sim} demonstrate that equivalent rates to those of \cite{zhan:16} are obtained for proper bandwidth choices that force the extra term in Theorem~\ref{thm: meanGen} to be of smaller order than the others.
\end{rem}

\begin{rem}
\label{rem: errAss}
Assumption C3 stipulates that the errors $\epsijl$ are iid and sub-Exponential across $i$ and $\ell$ for each fixed $j$, but dependence across $j$ is arbitrary.  In fact, the same rate obtained in Theorem~\ref{thm: meanGen} will hold under the weaker assumption $\epsijl = \epsilon_{ij}(\Tijl)$ for iid processes $\epsilon_{ij}$ satisfying $\max_{j \in \{1,\ldots,p\}} \sup_{t \in \Tj} E(\exp\{\lambda \epsilon_{ij}(t)\}) \leq \exp\{\lambda^2\sigma^2/2\}$ for some $0 < \sigma^2 < \infty$ and all $|\lambda| < (\sigma^2)\inv.$  The reason that this does not affect the rate is that the error dependence does not dominate than the intracurve dependence present in the latent functional data $X_i$. 
\end{rem}

Explicit rates will now be presented when a common weighting scheme, OBS or SUBJ, is used for all curves, although Theorem~\ref{thm: meanGen} can still be applied if different weighting schemes are applied to different indices $j$.  \cite{zhan:16} demonstrated that the relative performance of OBS and SUBJ depends on the values of the $\Nij,$ which may have different behaviors for different indices $j.$  However, in practice, it can be difficult to determine which weighting scheme is best for a given curve, except for an extreme case where some indices $j$ have much larger number of observations per curve compared to others.  Thus, for clarity, the weights $\wij$ are hereafter assumed to be constructed according to either OBS or SUBJ for all $j$.  Define $\Njtbar = n\inv \son \Nij^2$, $\Njp = \max_{i \in \{1,\ldots,n\}} \Nij,$ and $\NjH = n(\son \Nij\inv)\inv$.

\begin{cor}
\label{cor: meanObsSubjL2}
Suppose the assumptions of Theorem~\ref{thm: meanGen} hold.
\begin{enumerate}
    \item[{\upshape (i)}] OBS: 
    \[
    \max_{j \in \{1,\ldots,p\}} \normj{\hmujo - \muj} = O_P\left(\max_{j \in \{1,\ldots,p\}} \left[\bmuj^2 + \left\{\frac{\log(p)}{n}\left(\frac{1}{\Njbar \bmuj} + \frac{\Njtbar}{(\Njbar)^2}\right)\right\}^{1/2} + \frac{\Njp\log(p)}{n\Njbar \bmuj}\right] \right).
    \]
    \item[{\upshape (ii)}] SUBJ:
    \[
    \max_{j \in \{1,\ldots,p\}} \normj{\hmujs - \muj} = O_P\left(\max_{j \in \{1,\ldots,p\}} \left[\bmuj^2 + \left\{\frac{\log(p)}{n}\left(\frac{1}{\NjH \bmuj} + 1\right)\right\}^{1/2} + \frac{\log(p)}{n\bmuj}\right] \right).
    \]
\end{enumerate}
\end{cor}

\begin{rem}
\label{rem: compWt}
The rates of Corollary~\ref{cor: meanObsSubjL2} allow one to distinguish between settings in which each weighting scheme is expected to outperform the other, in line with the findings of \cite{zhan:16}.  Specifically, if the $N_{ij}$ are sufficiently homogeneous across $i$ for each $j$, in the sense that 
\begin{equation}
\label{eq: OBShom}
\lsn \max_{j \in \{1,\ldots,p\}} \max\left\{\frac{\Njtbar}{\Njbar^{2}}, \frac{\Njp}{\Njbar}\right\} < \infty,
\end{equation}
then the OBS scheme is never worse than the SUBJ scheme due to the fact that $\Njbar \geq \NjH.$  Indeed, it is possible that the above homogeneity condition holds and that, in addition, $\max_{j \in \{1,\ldots,p\}} \NjH = O(1)$ while $\min_{j \in \{1,\ldots,p\}} \Njbar \rightarrow \infty,$ in which case the OBS scheme is strictly better than SUBJ.  However, if such homogeneity fails, the OBS scheme rate can suffer considerably, as it can place too much weight on a small proportion of curves. An interesting suggestion of \cite{zhan:16} (see Remark 6 therein) for an alternative weighting scheme was to choose an optimal convex combination of the OBS and SUBJ weighting schemes, leading to an asymptotic rate that is better than either, although this theoretical appeal did not always lead to finite sample improvements in their numerical experiments.
\end{rem}

In the case of diverging $p$, Corollary~\ref{cor: meanObsSubjL2} suggests slightly inflated optimal bandwidth choices compared to the fixed $p$ case.  When the the curves for different indices $j$ are a mixture of non-dense, dense, and ultra-dense, the rate is driven by the worst rate of any individual mean curve estimate.  Under an FR design, an additional complication associated with the maximum rate over $j$ for diverging $p$ is that constants cannot be always ignored in the specification of optimal bandwidths.  Thus, for clarity, the following results provide optimal bandwidth choices in accordance with the behavior of the most sparsely and most densely observed functions, represented by $\Nmin = \min_{j \in \{1,\ldots,p\}} \Njbar$ and $\Nmax = \max_{j \in \{1,\ldots,p\}} \Njbar$ for the OBS scheme or $\NHmin = \min_{j \in \{1,\ldots,p\}} \NjH$ and $\NHmax = \max_{j \in \{1,\ldots,p\}} \NjH$ for the SUBJ scheme.  Remarkably, regardless of the type of functional data available, consistent estimation is possible whenever $\log(p)/n \rightarrow 0,$ so that exponential growth of $p$ is feasible for non-dense observations, including truly sparse ones in which the $\Nij$ are bounded.  For positive sequences $a_{jn}$ and $b_{jn},$ denote by $a_{jn} \preceq b_{jn}$ and $a_{jn} \asymplu b_{jn}$ the conditions $\lsn a_{jn}b_{jn}\inv < \infty$ and $\lsn \max_{j \in \{1,\ldots,p\}} a_{jn}b_{jn}\inv < \infty$, respectively. If $a_{jn} \preceq b_{jn}$ and $b_{jn} \preceq a_{jn}$ (respectively, $a_{jn} \asymplu b_{jn}$ and $b_{jn} \asymplu a_{jn}$), write $a_{jn} \asymp b_{jn}$ (resp., $a_{jn} \asympu b_{jn}$).

\begin{cor}
\label{cor: meanOptL2}
Assume an FR design and that assumptions A1, B1--B3, C1, and C2 hold.
\begin{enumerate}
    \item[{\upshape (i)}] OBS: Assume that \eqref{eq: OBShom} holds.
    \begin{enumerate}
        \item[{\upshape (a)}] If $\Nmax\left\{\frac{\log(p)}{n}\right\}^{1/4} \rightarrow 0$ and $\bmuj \asympu \left\{\frac{\log(p)}{n\Njbar}\right\}^{1/5}$, then
        $
        \max_{j \in \{1,\ldots,p\}} \normj{\hmujo - \muj} = O_P\left[\left\{\frac{\log(p)}{n\Nmin}\right\}^{2/5}\right].
        $
        \item[{\upshape (b)}] If $0 < \liminf_{n \rightarrow \infty} \Nmin\left\{\frac{\log(p)}{n}\right\}^{1/4} \leq \lsn \Nmax\left\{\frac{\log(p)}{n}\right\}^{1/4} < \infty$ and $\bmuj \asympu \left\{\frac{\log(p)}{n}\right\}^{1/4},$ then \newline
        $
        \max_{j \in \{1,\ldots,p\}} \normj{\hmujo - \muj} = O_P\left[\left\{\frac{\log(p)}{n}\right\}^{1/2}\right].
        $
        \item[{\upshape (c)}] If $\Nmin \left\{\frac{\log(p)}{n}\right\}^{1/4} \rightarrow \infty$, $\bmuj \asympu b_n,$ where $b_n = o\left[\left\{\frac{\log(p)}{n}\right\}^{1/4}\right]$, $b_n\Nmin \rightarrow \infty,$ and $b_n\left\{\frac{\log (p)}{n}\right\}^{-1/2} \rightarrow \infty,$ then
        $
        \max_{j \in \{1,\ldots,p\}} \normj{\hmujo - \muj} = O_P\left[\left\{\frac{\log(p)}{n}\right\}^{1/2}\right].
        $
    \end{enumerate}
    \item[{\upshape (ii)}]SUBJ: Replacing $\Njbar$, $\Nmin$, and $\Nmax$ with $\NjH$, $\NHmin$, and $\NHmax$, respectively, in parts {\upshape (a)--(c)} of {\upshape (i)}  leads to the corresponding results for $\hmujs$.
\end{enumerate}
\end{cor}

\begin{rem}
\label{rem: balance}
Compared to the optimal bandwidth choices outlined in \cite{zhan:16} for univariate functional data in the non-dense, dense, and ultra-dense regimes, those of Corollary~\ref{cor: meanOptL2} are nearly identical except that the division sample size $n$ is replaced by $n/\log(p)$ in the above result.  One minor exception is in the ultra-dense case (c); since the extra term $\log(p)/(n\bmuj)$ in Corollary~\ref{cor: meanObsSubjL2} is not affected by $\Njbar,$ under the FR design it is necessary to add an additional constraint to ensure that the smallest bandwidth does not decay more rapidly than $\{\log(p)/n\}^{1/2}$.  In doing so, this tail term is always of a smaller order than the others, so does not affect the final rate.  In the non-dense case (a), it is the bias and first stochastic term that dominate since $\Njbar \bmuj \rightarrow 0$, and the rate can be anywhere between $\left\{\log(p)/n\right\}^{2/5}$ (inclusive) and $\left\{\log(p)/n\right\}^{1/2}$ (exclusive).  In the dense case, the bias and stochastic terms are all balanced; only in the ultra-dense case do the terms involving the bandwidth become inconsequential, all converging at a rate faster than $\{\log(p)/n\}^{1/2}.$ 
\end{rem}

Lastly, under an SR design, similar rates are obtained under simpler assumptions.  Let $\Nbar$, $\NH$, $\Ntbar$, and $\Np$ be the common values of $\Njbar$, $\NjH$, $\Njtbar,$ and $\Njp$, respectively, across $j$.  
\begin{cor}
\label{cor: meanOptL2_sim}
Assume an SR design and that assumptions A1, B1--B3, C1, and C2 hold.
\begin{enumerate}
    \item[{\upshape (i)}] OBS: Assume that \eqref{eq: OBShom} holds.
    \begin{enumerate}
        \item[{\upshape (a)}] If $\Nbar\left\{\frac{\log(p)}{n}\right\}^{1/4} \rightarrow 0$, $\bmu \asymp \left\{\frac{\log(p)}{n\Nbar}\right\}^{1/5}$, then
        $
        \max_{j \in \{1,\ldots,p\}} \normj{\hmujo - \muj} = O_P\left[\left\{\frac{\log(p)}{n\Nbar}\right\}^{2/5}\right].
        $
        \item[{\upshape (b)}] If $\Nbar\left\{\frac{\log(p)}{n}\right\}^{1/4} \rightarrow C \in (0, \infty)$ and $\bmu \asymp \left\{\frac{\log(p)}{n}\right\}^{1/4},$ then
        $
        \max_{j \in \{1,\ldots,p\}} \normj{\hmujo - \muj} = O_P\left[\left\{\frac{\log(p)}{n}\right\}^{1/2}\right].
        $
        \item[{\upshape (c)}] If $\Nbar \left\{\frac{\log(p)}{n}\right\}^{1/4} \rightarrow \infty$, $\bmu = o\left[\left\{\frac{\log(p)}{n}\right\}^{1/4}\right]$ and $\bmu\Nbar \rightarrow \infty,$ then
        $
        \max_{j \in \{1,\ldots,p\}} \normj{\hmujo - \muj} = O_P\left[\left\{\frac{\log(p)}{n}\right\}^{1/2}\right].
        $
    \end{enumerate}
    \item[{\upshape (ii)}]SUBJ: Replacing $\Nbar$ with $\NH$ in parts {\upshape (a)--(c)} of {\upshape (i)} leads to the corresponding results for $\hmujs$.
\end{enumerate}
\end{cor}

\subsection{Uniform Convergence\label{ss: Supmean}}

Next, consider rates in the uniform norms $\normsupj{\cdot}.$  Following \cite{li:10} and \cite{zhan:16}, these rates will be given in the almost sure sense in order to allow for an easy comparison. 
 Remarks~\ref{rem: design}--\ref{rem: errAss} also apply to the following result, which provides the rate of convergence for a general weighting scheme.
\begin{thm}
\label{thm: meanSupGen}
Under assumptions A1, A2, B1--B3, and C1--C4, almost surely,
\[
\max_{j \in \{1,\ldots,p\}} \normsupj{\hmuj - \muj} = O\left(\max_{j \in \{1,\ldots,p\}} \left[\bmuj^2 + \left\{\log(p\vee n)\son\wij^2\Nij(\bmuj\inv + \Nij - 1)\right\}^{1/2} + \log(p \vee n)\bmuj\inv\wjn\right]\right).
\]
\end{thm}

\begin{rem}
\label{rem: supNorm}
    As mentioned briefly at the beginning of this section, Theorem~\ref{thm: meanSupGen} does not require any moment bounds on $\normsupj{X_{ij} - \muj}$, an assumption employed in \cite{li:10,zhan:16}; indeed, $X_{ij}$ need not even have bounded sample paths in order for the above result to hold.  In the proof, one approximates the supremum by a maximum over a grid, with the error in this approximation being bounded by the grid mesh size times a multiple of $\son \wij \soNij |U_{ij}(\Tijl)|.$  Under the sub-Exponential condition in assumption C1, this can be controlled directly without using the bound $|U_{ij}(\Tijl)| \leq \normsupj{U_{ij}},$ regardless of the strength of intracurve dependence present in the latent functional data.
\end{rem}

\begin{rem}
    \label{rem: pveen}
    The rates derived in \cite{li:10} and \cite{zhan:16} for $p = 1$ are nearly the same rates as those for weak consistency in the $L^2$ norm, with the only difference being a slight inflation by $\{\log(n)\}^{1/2}$.   Theorem~\ref{thm: meanSupGen} demonstrates a similar phenomenon in the high-dimensional setting, where the rates in Section~\ref{ss: L2mean} are augmented by $\{\log(p \vee n)\}$. Typically, the results derived in this paper will be of greatest interest when $p$ is much larger than $n$, in which case the almost sure uniform rates matches exactly the weak consistency rates in the $L^2$ norm. However, in case $n/p \rightarrow \infty,$ the new result demonstrates that the uniform rates are strictly worse than the $L^2$ rates, as should be expected.
\end{rem}

The corresponding rates for OBS and SUBJ weighting schemes and the corresponding divisions into non-dense, dense, and ultra-dense data are immediate.  Again, Remarks~\ref{rem: compWt} and \ref{rem: balance} also apply to the following corollaries.

\begin{cor}
\label{cor: meanObsSubjSup}
Suppose the assumptions of Theorem~\ref{thm: meanSupGen} hold.
\begin{enumerate}
    \item[{\upshape (i)}] OBS: Almost surely,
    \[
    \max_{j \in \{1,\ldots,p\}} \normsupj{\hmujo - \muj} = O\left\{\max_{j \in \{1,\ldots,p\}} \Bigg(\bmuj^2 + \left[\frac{\log(p\vee n)}{n}\left\{\frac{1}{\Njbar \bmuj} + \frac{\Njbar}{(\Njbar)^2}\right\}\right]^{1/2} + \frac{\Njp\log(p \vee n)}{n\Njbar \bmuj}\Bigg) \right\}.
    \]
    \item[{\upshape (ii)}] SUBJ: Almost surely,
    \[
    \max_{j \in \{1,\ldots,p\}} \normsupj{\hmujs - \muj} = O\left(\max_{j \in \{1,\ldots,p\}} \Bigg[\bmuj^2 + \left\{\frac{\log(p \vee n)}{n}\left(\frac{1}{\NjH \bmuj} + 1\right)\right\}^{1/2} + \frac{\log(p \vee n)}{n\bmuj}\Bigg] \right).
    \]
\end{enumerate}
\end{cor}

\begin{cor}
\label{cor: meanOptSup}
Assume an FR design and that assumptions A1, A2, B1--B3, C1, and C2 hold.
\begin{enumerate}
    \item[{\upshape (i)}] OBS: Assume that \eqref{eq: OBShom} holds.
    \begin{enumerate}
        \item[{\upshape (a)}] If $\Nmax\left\{\frac{\log(p \vee n)}{n}\right\}^{1/4} \rightarrow 0$ and $\bmuj \asympu \left\{\frac{\log(p\vee n)}{n\Njbar}\right\}^{1/5}$, then
        $
        \max_{j \in \{1,\ldots,p\}} \normsupj{\hmujo - \muj} = O\left[\left\{\frac{\log(p \vee n)}{n\Nmin}\right\}^{2/5}\right]
        $ \\ almost surely.
        \item[{\upshape (b)}] If $0 < \liminf_{n\rightarrow \infty} \Nmin\left\{\frac{\log(p \vee n)}{n}\right\}^{1/4} \leq \lsn \Nmax\left\{\frac{\log(p \vee n)}{n}\right\}^{1/4} < \infty$ and $\bmuj \asympu \left\{\frac{\log(p \vee n)}{n}\right\}^{1/4},$ \\ then
        $
        \max_{j \in \{1,\ldots,p\}} \normsupj{\hmujo - \muj} = O\left[\left\{\frac{\log(p\vee n)}{n}\right\}^{1/2}\right]
        $ almost surely.
        \item[{\upshape (c)}] If $\Nmin \left\{\frac{\log(p \vee n)}{n}\right\}^{1/4} \rightarrow \infty$, $\bmuj \asympu b_n$, where $b_n = o\left[\left\{\frac{\log(p\vee n)}{n}\right\}^{1/4}\right]$, $b_n\Nmin \rightarrow \infty,$ and \\ $b_n\left\{\frac{\log(p \vee n)}{n}\right\}^{-1/2}\rightarrow \infty,$  then
        $
        \max_{j \in \{1,\ldots,p\}} \normsupj{\hmujo - \muj} = O\left[\left\{\frac{\log(p \vee n)}{n}\right\}^{1/2}\right]
        $ almost surely.
    \end{enumerate}
    \item[{\upshape (ii)}]SUBJ: Replacing $\Njbar$, $\Nmin$, and $\Nmax$ with $\NjH$, $\NHmin$, and $\NHmax$, respectively, in parts {\upshape (a)--(c)} of {\upshape (i)}  leads to the corresponding results for $\hmujs$.
\end{enumerate}
\end{cor}

\begin{cor}
\label{cor: meanOptSupSim}
Assume an SR design and that assumptions A1, A2, B1--B3, C1, and C2 hold.
\begin{enumerate}
    \item[{\upshape (i)}] OBS: Assume that \eqref{eq: OBShom} holds.
    \begin{enumerate}
        \item[{\upshape (a)}] If $\Nbar\left\{\frac{\log(p \vee n)}{n}\right\}^{1/4} \rightarrow 0$ and $\bmu \asymp \left\{\frac{\log(p \vee n)}{n\Nbar}\right\}^{-1/5}$, then
        $
        \max_{j \in \{1,\ldots,p\}} \normsupj{\hmujo - \muj} = O\left[\left\{\frac{\log(p \vee n)}{n\Nbar}\right\}^{2/5}\right]
        $ \\ almost surely.
        \item[{\upshape (b)}] If $\Nbar\left\{\frac{\log(p \vee n)}{n}\right\}^{1/4} \rightarrow C \in (0, \infty)$ and $\bmu \asymp \left\{\frac{\log(p \vee n)}{n}\right\}^{1/4},$ then, almost surely, \\
        $
        \max_{j \in \{1,\ldots,p\}} \normsupj{\hmujo - \muj} = O\left[\left\{\frac{\log(p\vee n)}{n}\right\}^{1/2}\right].$
        \item[{\upshape (c)}] If $\Nbar \left\{\frac{\log(p \vee n)}{n}\right\}^{1/4} \rightarrow \infty$, $\bmu = o\left[\left\{\frac{\log(p \vee n)}{n}\right\}^{1/4}\right]$ and $\bmu\Nbar \rightarrow \infty,$ then, almost surely, \\
        $
        \max_{j \in \{1,\ldots,p\}} \normsupj{\hmujo - \muj} = O\left[\left\{\frac{\log(p \vee n)}{n}\right\}^{1/2}\right]
        $.
    \end{enumerate}
    \item[{\upshape (ii)}]SUBJ: Replacing $\Nbar$ with $\NH$ in parts {\upshape (a)--(c)} of {\upshape (i)} leads to the corresponding results for $\hmujs$.
\end{enumerate}
\end{cor}

\section{Covariance Estimation\label{sec: cov}}

Next, rates of convergence for the covariance estimators will be provided.  In order to separate the effects of mean estimation from covariance estimation, the mean is assumed to be known so that the raw covariances $Z_{ijklm}$ used in \eqref{eq: CovEst} are replaced by $\Zijklm = (\Yijl - \muj(\Tijl))(\Yikm - \muk(\Tikm)).$  The ensuing estimates will still be denoted $\hgammajk$.  The rates for the true empirical estimator are obtained by adding the uniform mean convergence rates, since raw covariances involve pointwise evaluations of the mean estimates.

\subsection{\texorpdfstring{$L^2$}{L2} Convergence\label{ss: covL2}}

Under an FR design, define $\vjkn = \max_{i \in \{1,\ldots,n\}} \vijk|\Iijk|$ and
\begin{equation}
\label{eq: qjkn}
\qjkn = \begin{cases}
    \son \vijk^2\Nij\Nik(\bgammaj\inv + \Nij - 1)(\bgammak\inv + \Nik - 1), & j \neq k, \\
    \son \vijj^2\Nij(\Nij - 1)\left\{\bgammaj^{-2} + 2\bgammaj\inv(\Nij - 2) +(\Nij - 2)(\Nij - 3)\right\}, & j = k.
\end{cases}
\end{equation}
Under an SR design, for all $j,k$, take $\vjkn = \max_{i \in \{1,\ldots,n\}} v_i|\mathcal{I}_i|$ and $\qjkn$ to be the common value of $q_{jjn}$ across $j$ in \eqref{eq: qjkn}.  In addition, define
\begin{equation}
    \label{eq: omegajk}
    \omega_{jk}^2 = \begin{cases}
        \max(\son \vijk^2\Nij^2 \Nik, \son \vijk^2 \Nij \Nik^2), & j \neq k, \\
        \son \vijj^2\Nij(\Nij - 1)^2, & j = k,
    \end{cases}
\end{equation}
and let $\omega$ denote the common value of $\omega_{jj}$ under an SR design.  The following assumptions will be used.
\begin{enumerate}
    \item[D1] For each $j \in \{1,\ldots,p\}$ and $t \in \Tj,$ $\Uij(t)$ is sub-Gaussian with parameter $\theta_j(t) > 0$, that is, $E[\exp\{\lambda \Uij(t)\}] \leq \exp\{\lambda^2\theta_j^2(t)/2\}$ for all $\lambda \in \mathbb{R}.$  Furthermore, 
    $
    \overline{\theta} = \lim_{n \rightarrow \infty} \max_{j \in \{1,\ldots,p\}} \sup_{t \in \Tj} \theta_j(t) < \infty.
    $
    \item[D2] For each $j$, the errors $\epsijl$ are sub-Gaussian random variables with parameter $\sigma_j,$ that is, $E[\exp\{\lambda \epsijl\}] \leq \exp\{\lambda^2\sigma_j^2/2\}$ for all $\lambda \in \mathbb{R}.$  Furthermore,
    $
    \sigma = \lim_{n\rightarrow \infty} \max_{j \in \{1,\ldots,p\}} \sigma_j < \infty.
    $
    In addition, under an FR design, the $\epsijl$ are independent across all indices; under an SR design, $\epsijl$ and $\epsilon_{i'j'\ell'}$ are independent whenever $(i,\ell) \neq (i',\ell')$.
    \item[D3] The bandwidths satisfy $\max_{j \in \{1,\ldots,p\}} \bgammaj \rightarrow 0$ and $\log(p)\max_{j,k \in \{1,\ldots,p\}} \qjkn \rightarrow 0$, with $\qjkn$ defined in \eqref{eq: qjkn}. In addition,  $\log(p)\max_{j,k \in \{1,\ldots,p\}} \bgammaj\inv\bgammak\inv\omega_{jk}^2\rightarrow 0$ under an FR design; under an SR design, $\bgamma^{-2} \omega^2 \rightarrow 0$.
    \item[D4] There exists $\alpha > 0$ such that, with $\qjkn$ as defined in \eqref{eq: qjkn},
    $$
    \frac{n^{\alpha}\max_{j,k \in \{1,\ldots,p\}} \left[\left\{\log(p\vee n)\qjkn\right\}^{1/2} + \log(p\vee n)\bgammaj\inv\bgammak\inv\vjkn\right]}{\max_{j \in \{1,\ldots,p\}} \bgammaj^{-3}} \rightarrow \infty.
    $$
    \end{enumerate}

Assumptions D1--D4 correspond to assumptions C1--C4, but adapted for covariance estimation.  Critically, assumptions D1 and D2 impose sub-Gaussian rather than sub-Exponential tails, due to the fact that covariance estimates involve averages of products, and products of sub-Gaussian random variables have sub-Exponential tails. Proofs of all results in this section can be found in Section~\ref{AppE}, with auxiliary lemmas and their proofs given in Section~\ref{AppD}.

\begin{thm}
\label{thm: CovGen}
Under assumptions A1, B1, B2, B4, and D1--D3, 
\[
\max_{j,k \in \{1,\ldots,p\}} \normjk{\hgammajk - \gammajk} = O_P\left(\max_{j \in \{1,\ldots,p\}} \bgammaj^2 + \max_{j,k \in \{1,\ldots,p\}} \left[\left\{\log(p)\qjkn\right\}^{1/2} + \log(p)\bgammaj\inv\bgammak\inv\vjkn \right]\right).
\]
\end{thm}

\begin{rem}
    \label{rem: designCov}
    As in Theorem~\ref{thm: meanGen}, the rates in Theorem~\ref{thm: CovGen} are stated in the context of an FR design, but also hold under an SR design using the common values of the relevant quantities across $j,k$, including a common bandwidth $\bgamma$, that have been defined previously.  Assumption D3 and comparison of Corollaries~\ref{cor: CovOptL2} and \ref{cor: CovOptL2_sim} below illustrate the weaker bandwidth requirements under an SR design. 
\end{rem}

\begin{rem}
    \label{rem: rateCompCov}
    Similar to Remark~\ref{rem: rateComp}, the rate in Theorem~\ref{thm: CovGen} can be compared with that of \cite{zhan:16}. As the latter did not consider cross-covariance estimation, consider $j = k.$  Once again, the arguments in the proof of Theorem~\ref{thm: CovGen}, when applied to a single auto-covariance estimate $\hat{\gamma}_{jj},$ lead to the rate $\bgammaj^2 + \qjkn^{1/2} + \bgammaj^{-2}\overline{v}_{njj},$ the first two terms matching the rate of \cite{zhan:16}.  As for mean estimation, the additional term in the derived rate arises from the involvement of higher order moments needed for the exponential tail bound.  This extra term has more impact on optimal bandwidth choice in covariance estimation than mean estimation, 
    but still does not affect the overall rates.
\end{rem}

The following corollaries translate the rate of Theorem~\ref{thm: CovGen} to the OBS and SUBJ weighting schemes under FR and SR designs, respectively.  For $r,s \in \{1,2\}$, define $\overline{N}_{j^rk^{s}} = n\inv\son \Nij^r\Nik^{s},$ $\Njkbar = \overline{N}_{j^1k^1},$ $\Njkp = \max_{i \in \{1,\ldots,n\}} \Nij\Nik$, $N_{j^rk^{s}}^H = n\left(\son \Nij^{-r}\Nik^{-s}\right)\inv,$ and $\NjkH = N_{j^1k^1}^H.$  Under an SR design, let $N_{(2)}^H$ denote the common value of $N_{jk}^H$ across $j,k$.

\begin{cor}
    \label{cor: CovObsSubjL2}
    Suppose the assumptions of Theorem~\ref{thm: CovGen} hold.  
    \begin{enumerate}
        \item[{\upshape (i)}] OBS:
        \[
        \begin{split}
        \max_{j,k \in \{1,\ldots,p\}} \normjk{\hgammajko - \gammajk} &= O_P\left(\max_{j \in \{1,\ldots,p\}} \bgammaj^2 + \max_{j,k \in \{1,\ldots,p\}}\left[\left\{\frac{\log(p)}{n}\left(\frac{1}{\Njkbar\bgammaj\bgammak} + \frac{\Njkkbar}{\Njkbar^2\bgammaj} +  
        \frac{\Njjkkbar}{\Njkbar^2}\right)\right\} ^{1/2} \right. \right. \\
        &\hspace{4.5cm} + \left. \left. \frac{\Njkp\log(p)}{n\Njkbar\bgammaj\bgammak}\right]\right).
        \end{split}
        \]
        \item[{\upshape (ii)}] SUBJ:
        \[
        \max_{j,k \in \{1,\ldots,p\}} \normjk{\hgammajko - \gammajk} = O_P\left(\max_{j \in \{1,\ldots,p\}} \bgammaj^2 + \max_{j,k \in \{1,\ldots,p\}}\left[\left\{\frac{\log(p)}{n}\left(\frac{1}{\NjkH\bgammaj\bgammak} + \frac{1}{\NjH\bgammaj} + 1\right)\right\}^{1/2} + \frac{\log(p)}{n\bgammaj\bgammak}\right]\right).
        \]
    \end{enumerate}
\end{cor}

\begin{rem}
\label{rem: compWtCov}
Corollary~\ref{cor: CovObsSubjL2} allows one to compare the two specific weighting schemes.  By strengthening the homogeneity condition on the observation numbers $\Nij$ in \eqref{eq: OBShom} to the condition
\begin{equation}
    \label{eq: OBShom2}
    \limsup_{n\rightarrow \infty} \max_{j,k \in \{1,\ldots,p\}}\max\left\{
    \frac{\Njbar\Njkkbar}{\Njkbar^2}, \frac{\Njjkkbar}{\Njkbar^2}, \frac{\Njkp}{\Njkbar}\right\} < \infty,
\end{equation}
then the OBS scheme is never worse than the SUBJ scheme due to the fact that $\Njbar \geq \NjH$ and $\Njkbar \geq \NjkH.$  If \eqref{eq: OBShom2} fails, however, the estimators may not be consistent under the OBS scheme. 
\end{rem}

\begin{cor}
\label{cor: CovOptL2}
Assume an FR design and that assumptions A1, B1, B2, B4, D1, and D2 hold.
\begin{enumerate}
    \item[{\upshape (i)}] OBS: Assume that \eqref{eq: OBShom2} holds and that $\limsup_{n\rightarrow \infty} \max_{j,k \in \{1,\ldots,p\}} (\Njbar\Nkbar)/\Njkbar < \infty.$
    \begin{enumerate}
        \item[{\upshape (a)}] If $\Nmax\left\{\frac{\log(p)}{n}\right\}^{1/4} \rightarrow 0$, $\bgammaj \asympu \left\{\frac{\log(p)}{n\Njbar^2}\right\}^{1/6},$ then
        $
        \max_{j,k \in \{1,\ldots,p\}} \normjk{\hgammajko - \gammajk} = O_P\left[\left\{\frac{\log(p)}{n\Nmin^2}\right\}^{1/3}\right].
        $
        \item[{\upshape (b)}] If $\liminf_{n \rightarrow \infty} \Nmin\left\{\frac{\log(p)}{n}\right\}^{1/4} > 0$ and $\bgammaj \asympu \left\{\frac{\log(p)}{n}\right\}^{1/4},$ then
        $
        \max_{j,k \in \{1,\ldots,p\}} \normjk{\hgammajko - \gammajk} = O_P\left[\left\{\frac{\log(p)}{n}\right\}^{1/2}\right].
        $
    \end{enumerate}
    \item[{\upshape (ii)}]SUBJ: Assume that $\limsup_{n\rightarrow \infty} \max_{j,k \in \{1,\ldots,p\}} (\NjH\NkH)/\NjkH < \infty.$ Replacing $\Njbar$, $\Nmin$, and $\Nmax$ with $\NjH$, $\NHmin$, and $\NHmax$, respectively, in parts {\upshape (a)} and {\upshape (b)} of {\upshape (i)} leads to the corresponding results for $\hgammajks$.
\end{enumerate}
\end{cor}

\begin{cor}
\label{cor: CovOptL2_sim}
Assume an SR design and that assumptions A1, B1, B2, B4, D1, and D2 hold.
\begin{enumerate}
    \item[{\upshape (i)}] OBS: Assume that \eqref{eq: OBShom2} holds.
    \begin{enumerate}
        \item[{\upshape (a)}] If $\Nbar\left\{\frac{\log(p)}{n}\right\}^{1/4} \rightarrow 0$, $\bgamma \asymp \left\{\frac{\log(p)}{n\Nbar^2}\right\}^{1/6}$, then
        $
        \max_{j,k \in \{1,\ldots,p\}} \normjk{\hgammajko - \gammajk} = O_P\left[\left\{\frac{\log(p)}{n\Nbar^2}\right\}^{1/3}\right].
        $
        \item[{\upshape (b)}] If $\liminf_{n \rightarrow \infty} \Nbar\left\{\frac{\log(p)}{n}\right\}^{1/4} > 0$ and $\bgamma \asymp \left\{\frac{\log(p)}{n}\right\}^{1/4},$ then
        $
        \max_{j,k \in \{1,\ldots,p\}} \normjk{\hgammajko - \gammajk} = O_P\left[\left\{\frac{\log(p)}{n}\right\}^{1/2}\right].
        $
    \end{enumerate}
    \item[{\upshape (ii)}]SUBJ: Assume that $\limsup_{n \rightarrow \infty} (\NH)^2/N_{(2)}^H < \infty$. Replacing $\Nbar$ with $\NH$ in parts {\upshape (a)} and {\upshape (b)} of {\upshape (i)} leads to the corresponding results for $\hgammajks$.
\end{enumerate}
\end{cor}

\begin{rem}
\label{rem: CovUD}
    Unlike Corollaries~\ref{cor: meanOptL2} and \ref{cor: meanOptL2_sim}, Corollaries~\ref{cor: CovOptL2} and \ref{cor: CovOptL2_sim} contain only two cases for each of the OBS and SUBJ weighting schemes, with the dense and ultra-dense regimes being combined.  This is due to the extra term in Corollary~\ref{cor: CovObsSubjL2} discussed previously in Remark~\ref{rem: rateCompCov}.  For mean estimation, in the ultra-dense regime, the bandwidths $\bmuj$ are all allowed to decay more quickly than $\{\log(p)/n\}^{1/4}.$  However, allowing the same behavior for covariance bandwidths $\bgammaj$ would cause the last term in the rate to be slower than $\{\log(p)/n\}^{1/2}$. Hence, for high-dimensional functional data, while this result does not distinguish between dense and ultra-dense observation designs in terms of covariance estimation, in both regimes one is still able to obtain the appropriate rate.  In the sparse regime, the rate can be anywhere between $\left\{\log(p)/n\right\}^{1/3}$ (inclusive) and $\left\{\log(p)/n\right\}^{1/2}$ (exclusive).
\end{rem}

\subsection{Uniform Convergence\label{ss: covUnif}}

For brevity, the strong uniform rates of convergence for covariance estimation will be stated without further discussion; the reader is referred to Remarks~\ref{rem: supNorm}, \ref{rem: pveen}, and \ref{rem: rateCompCov}--\ref{rem: CovUD} for relevant comments on these results.

\begin{thm}
\label{thm: CovSupGen}
Under assumptions A1, A2, B1, B2, B4, and D1--D4, almost surely,
\[
\max_{j,k \in \{1,\ldots,p\}} \normsupjk{\hgammajk - \gammajk}  = O\left(\max_{j \in \{1,\ldots,p\}} \bgammaj^2 + \max_{j,k \in \{1,\ldots,p\}} \left[\left\{\log(p\vee n)\qjkn\right\}^{1/2} + \log(p \vee n)\bgammaj\inv\bgammak\inv\vjkn \right]\right).
\]
\end{thm}

\begin{cor}
    \label{cor: CovObsSubjSup}
    Suppose the assumptions of Theorem~\ref{thm: CovGen} hold.  
    \begin{enumerate}
        \item[{\upshape (i)}] OBS:
        \[
        \begin{split}
        \max_{j,k \in \{1,\ldots,p\}} \normsupjk{\hgammajko - \gammajk} &= O_P\left(\max_{j \in \{1,\ldots,p\}} \bgammaj^2 + \max_{j,k \in \{1,\ldots,p\}}\left[\left\{\frac{\log(p \vee n)}{n}\left(\frac{1}{\Njkbar\bgammaj\bgammak} + \frac{\Njkkbar}{\Njkbar^2\bgammaj} + 
        \frac{\Njjkkbar}{\Njkbar^2}\right)\right\}^{1/2} \right. \right. \\ &\hspace{4.5cm} \left. \left. + \frac{\Njkp\log(p\vee n)}{n\Njkbar\bgammaj\bgammak}\right]\right).
        \end{split}
        \]
        \item[{\upshape (ii)}] SUBJ:
        \[
        \begin{split}
        \max_{j,k \in \{1,\ldots,p\}} \normsupjk{\hgammajko - \gammajk} &= O_P\left(\max_{j \in \{1,\ldots,p\}} \bgammaj^2 + \max_{j,k \in \{1,\ldots,p\}}\left[\left\{\frac{\log(p \vee n)}{n}\left(\frac{1}{\NjkH\bgammaj\bgammak} + \frac{1}{\NjH\bgammaj} + 1\right)\right\}^{1/2} \right. \right. \\ &\hspace{4.5cm} \left. \left. + \frac{\log(p \vee n)}{n\bgammaj\bgammak}\right]\right).
        \end{split}
        \]
    \end{enumerate}
\end{cor}

\begin{cor}
\label{cor: CovOptSup}
Assume an FR design and that assumptions A1, B1, B2, B4, D1, and D2 hold.
\begin{enumerate}
    \item[{\upshape (i)}] OBS: Assume that \eqref{eq: OBShom2} holds and that $\limsup_{n\rightarrow \infty} \max_{j,k \in \{1,\ldots,p\}} (\Njbar\Nkbar)/\Njkbar < \infty.$
    \begin{enumerate}
        \item[{\upshape (a)}] If $\Nmax\left\{\frac{\log(p\vee n)}{n}\right\}^{1/4} \rightarrow 0$, $\bgammaj \asympu \left\{\frac{\log(p\vee n)}{n\Njbar^2}\right\}^{1/6}$, then
        $
        \max_{j,k \in \{1,\ldots,p\}} \normsupjk{\hgammajko - \gammajk} = O\left[\left\{\frac{\log(p\vee n)}{n\Nmin^2}\right\}^{1/3}\right]
        $ \\ almost surely.
        \item[{\upshape (b)}] If $\liminf_{n \rightarrow \infty} \Nmin\left\{\frac{\log(p\vee n)}{n}\right\}^{1/4} >0 $ and $\bgammaj \asympu \left\{\frac{\log(p\vee n)}{n}\right\}^{1/4},$ then, almost surely, \\
        $
        \max_{j,k \in \{1,\ldots,p\}} \normsupjk{\hgammajko - \gammajk} = O\left[\left\{\frac{\log(p\vee n)}{n}\right\}^{1/2}\right]
        $.
    \end{enumerate}
    \item[{\upshape (ii)}]SUBJ: Assume that $\limsup_{n\rightarrow \infty} \max_{j,k \in \{1,\ldots,p\}} (\NjH\NkH)/\NjkH < \infty.$ Replacing $\Njbar$, $\Nmin$, and $\Nmax$ with $\NjH$, $\NHmin$, and $\NHmax$, respectively, in parts {\upshape (a)} and {\upshape (b)} of {\upshape (i)} leads to the corresponding results for $\hgammajks$.
\end{enumerate}
\end{cor}

\begin{cor}
\label{cor: CovOptSup_sim}
Assume an SR design and that assumptions A1, B1, B2, B4, D1, and D2 hold.
\begin{enumerate}
    \item[{\upshape (i)}] OBS: Assume that \eqref{eq: OBShom2} holds.
    \begin{enumerate}
        \item[{\upshape (a)}] If $\Nbar\left\{\frac{\log(p\vee n)}{n}\right\}^{1/4} \rightarrow 0$, and $\bgamma \asymp \left\{\frac{\log(p\vee n)}{n\Nbar^2}\right\}^{1/6}$, then
        $
        \max_{j,k \in \{1,\ldots,p\}} \normjk{\hgammajko - \gammajk} = O\left[\left\{\frac{\log(p\vee n)}{n\Nbar^2}\right\}^{1/3}\right]
        $ \\ almost surely.
        \item[{\upshape (b)}] If $\liminf_{n \rightarrow \infty} \Nbar\left\{\frac{\log(p\vee n)}{n}\right\}^{1/4} > 0$ and $\bgamma \asymp \left\{\frac{\log(p\vee n)}{n}\right\}^{1/4},$ then, almost surely, \\
        $
        \max_{j,k \in \{1,\ldots,p\}} \normjk{\hgammajko - \gammajk} = O\left[\left\{\frac{\log(p\vee n)}{n}\right\}^{1/2}\right]
        $.
    \end{enumerate}
    \item[{\upshape (ii)}]SUBJ: Assume that $\limsup_{n \rightarrow \infty} (\NH)^2/N_{(2)}^H < \infty$. Replacing $\Nbar$ with $\NH$ in parts {\upshape (a)} and {\upshape (b)} of {\upshape (i)} leads to the corresponding results for $\hgammajks$.
\end{enumerate}
\end{cor}

\section{Discussion\label{sec: disc}}

The results derived in this paper provide an important foundation for high-dimensional functional data analysis by establishing sufficient conditions under which uniform consistency of a diverging number of mean and covariance estimates will hold.  Importantly, the practical reality of discrete and noisy functional observations does not preclude consistent estimation in high dimensions; indeed, the results lead to a natural division of high-dimensional functional data into three regimes (sparse, dense, and ultra-dense) based on the behavior of the average number of observations available per component curve relative to $\{\log(p)/n\}^{1/4},$ providing the expected generalization of the regime divisions discovered in \cite{zhan:16} for univariate functional data.  By properly choosing the smoothing bandwidths, the worst case scenario for sparsely observed functions is the optimal nonparametric high-dimensional convergence rate for both mean and covariance estimation; for densely or ultra-densely observed curves, the parametric rate is always attainable.  

The results utilize concentration inequalities in Hilbert spaces \citep{bosq:00} that require stricter tail assumptions on the pointwise behavior of the functional data compared to previous results for the case $p = 1$.  In this paper, sub-Exponential and sub-Gaussian tails were used for mean and covariance estimation, respectively.  With these tail assumptions, the results were able to successfully distinguish between dense and ultra-dense functional data for both the OBS or SUBJ weighting scheme in mean estimation, but not so for covariance estimation, due to an additional term in the rate that was not present in previous work \citep{li:10,zhan:16} in the case $p = 1$; see Remarks~\ref{rem: rateCompCov} and \ref{rem: CovUD}.  Thus, it is possible that stronger concentration inequalities or alternative modes of analysis may lead to such a distinction in future work. For example, under the assumption of tails with only polynomial decay, similar arguments to those of \cite{zhan:16} could be used to derive rates that, while slower than those obtained here, have the potential to effectively distinguish between dense and ultra-dense data.  Besides different tail assumptions, there are many other choices that may lead to different divisions, including different learning tasks such as regression or classification, modes of convergence (e.g., pointwise asymptotic normality), or quantification of estimation error besides the $L^2$ and uniform metrics.  Moreover, while the theorems and corrollaries present asymptotic properties, the lemmas in Sections~\ref{AppB} and \ref{AppD} are non-asymptotic in nature, and may be useful for finite sample inference.

The use of local linear methodology in the analysis of high-dimensional functional data analysis will undoubtedly have many limitations, some of which are theoretical and others that are practical.  This paper has addressed the former via asymptotic analyses as well as via simulations in the supplementary material. Important practical questions, such as bandwidth selection or approximations for computational speedup, will require adaptation of the usual tools for univariate functional data if they are to remain scalable.  For a brief discussion of some of these issues, see Section 1.2 in the supplementary material.

Although the terms sparse and dense appear frequently in this paper, there has been no application of regularization or shrinkage in the estimation procedure.  It is well-known that, for multivariate non-functional data of high dimension, such approaches can lead to improved rates of convergence and enhanced interpretability.  For high-dimensional functional data, such regularization has been successfully applied in regression models and graphical model estimation, although theoretical guarantees have mostly been established in the oracle case of fully observed functions. The results in this paper will provide a path for theoretical investigation of these and other novel methodologies for high-dimensional functional data that are applicable under any observational design, requiring only initial mean and covariance plug-in estimates.

\section{Technical Details\label{sec: App}}

This section provides technical arguments for the theoretical results stated in Sections~\ref{sec: mean} and \ref{sec: cov}.  

\subsection{Auxiliary Lemmas for Mean Estimation\label{AppB}}

\input{MeanLemmasJMVA.tex}

\subsection{Proofs of Results in Section~\ref{sec: mean}\label{AppC}}

\input{MeanTheoremsAndCorollariesJMVA.tex}

\subsection{Auxiliary Lemmas For Covariance Estimation\label{AppD}}

\input{CovLemmasJMVA.tex}

\subsection{Proofs of Results in Section~\ref{sec: cov}\label{AppE}}

\input{CovTheoremsAndCorollariesJMVA.tex}

\section*{Acknowledgments}

This research was supported by NSF Award DMS-2310943.

\section*{Supplementary Materials} 

The supplementary materials include comparisons of the performance of the OBS and SUBJ schemes on simulated data sets, along with corresponding code.

\bibliographystyle{myjmva}


\end{document}

%% file: MeanLemmasJMVA.tex
For $r\in\{0,1,2\},$ and $j\in \{1,\dots, p\},$ define
\begin{equation*}
    S_{jr}(t) = \son \wij \soNij \Kb{\bmuj}\pr{\Tijl-t}\pr{\frac{\Tijl-t}{\bmuj}}^r, \\ \quad
    R_{jr}(t) = \son \wij \soNij \Kb{\bmuj}\pr{\Tijl-t}\pr{\frac{\Tijl-t}{\bmuj}}^r \Yijl.
\end{equation*}
Then the error in mean estimation can be expressed as
\begin{equation}
    \label{eq: meanDiff}
    \hmuj(t) -\muj(t) = \frac{\left\{\Rzj(t) -\muj(t)\Szj(t) - \bmuj\muj'(t)\Soj(t)\right\}\Stj(t)}{\Szj(t)\Stj(t) - \Soj^2(t)}  - \frac{\left\{\Roj(t) - \muj(t)\Soj(t) - \bmuj\muj'(t)\Stj(t)\right\}\Soj(t)}{\Szj(t)\Stj(t) - \Soj^2(t)}.
\end{equation}
With $\Uijl = \Uij(\Tijl) + \epsijl$, the numerator in the first term on the right-hand side of \eqref{eq: meanDiff} can be expressed as
    \begin{equation}
    \label{eq: MnNumTermExp1}
    \begin{split}
    \Rzj(t) -\muj(t)\Szj(t) - \bmuj\muj'(t)\Soj(t) &=  \son \wij \soNij \Kb{\bmuj}(\Tijl - t)\Uijl \\
    & \hspace{0.6cm} + \son \wij \soNij \Kb{\bmuj}(\Tijl - t)(\Tijl - t)^2\int_0^1 \muj''\{t + v(\Tijl - t)\}(1 - v)\dv,
    \end{split}
    \end{equation}
which follows from a Taylor expansion.  Similarly,
    \begin{equation}
    \label{eq: MnNumTermExp2}
    \begin{split}
    \Roj(t) -\muj(t)\Soj(t) - \bmuj\muj'(t)\Stj(t) &= \son \wij \soNij \Kb{\bmuj}(\Tijl - t)\left(\frac{\Tijl - t}{\bmuj}\right)\Uijl \\
    & \hspace{0.6cm} + \son \wij \soNij \Kb{\bmuj}(\Tijl - t)\left(\frac{\Tijl - t}{\bmuj}\right)(\Tijl - t)^2\int_0^1 \muj''(t + v(\Tijl - t))(1 - v)\dv.
    \end{split}
    \end{equation}
    
\subsubsection{Uniform Convergence of \texorpdfstring{$S_{jr}$}{Sjr}, \texorpdfstring{$r \in \{0,1,2\}$}{r = 0,1,2}}
The results of this section allows for the derivation of lower bounds for the denominator of \eqref{eq: meanDiff} that hold uniformly in $t$ and across $j$ with high probability.  In this result and those that follow, constants that are uniform over all data generating mechanisms satisfying the mentioned assumptions will be denoted by $C_a$, $a \in \mathbb{N}$, and may take different values across different inequalities.  Finally, these results do not depend on whether the observation times follow a fully or simultaneous random design.

\begin{lma}
\label{lma: SjrConcIneq}
Suppose that assumptions A1 and B1 hold.  There exist $C_1, C_2 > 0$ such that, for any $\epsilon > 0$,
$$
\Pr\left(\normsupj{S_{jr}(\cdot) - E[S_{jr}(\cdot)]} > \epsilon\right) \leq C_1\exp\left\{\frac{-C_2\bmuj^2\epsilon^2}{\sum_{i = 1}^nw_{ij}^2 N_{ij}}\right\}, \quad r \in \{0, 1, 2\},\, j \in \{1,\ldots,p\}. 
$$
\end{lma}

\begin{proof}[\textbf{\upshape Proof:}]
    Define the weighted empirical distribution function $\hat{F}_j(t) = \son \wij \soNij \mathbf{1}(\Tijl \leq t).$ Take $V_{K,r}$ to be the total variation of the kernel $K(u)u^r$ over $[-1,1]$ for $r \in \{0,1,2\},$ an set $V_K = \max_{r \in \{0,1,2\}} V_{K,r}.$  Then standard arguments show that
    $
    \normsupj{S_{jr}(\cdot) - E[S_{jr}(\cdot)]} \leq V_K\bmuj\inv\normsupj{\hat{F}_j - F_j}.
    $
    Next, Lemma~1.1 of \cite{marc:84} implies that, for any $\tau > 0,$
    $$
    \Pr\left\{\normsupj{\hat{F}_j - F_j} > \tau\left(\son \wij^2 \Nij\right)^{1/2}\right\} \leq \left(1 + 2\sqrt{2\pi}\tau\right)\exp\left\{-\tau^2/8\right\}. 
    $$
    Since, for any polynomial $p(u),$ $\sup_{u > 0} p(u)e^{-u^2} \leq a_1e^{-a_2u^2}$ for some constants $a_1,$ $a_2 > 0$, it follows that, for any $\epsilon > 0,$ there are constants $C_1, C_2 > 0$ such that, as claimed,
    \[
        \Pr\left(\normsupj{S_{jr}(\cdot) - E[S_{jr}(\cdot)]} > \epsilon\right) \leq \Pr\left(\normsupj{\hat{F}_j - F_j} > V_K\inv \bmuj \epsilon\right)
        \leq C_1\exp\left\{\frac{-C_2\bmuj^2\epsilon^2}{\son \wij^2\Nij}\right\}.
    \]
\end{proof}

\begin{lma}
\label{lma: InfDen}
    Suppose assumptions A1, B1, B2, and C3 hold.  Then there exist $C_1, C_2, \eta, N > 0$ such that, for any $n \geq N$,
    \begin{equation*}
        \Pr\left[\inf_{t \in \Tj}\left\{\Szj(t)\Stj(t) - \Soj^2(t)\right\} < \eta\right] \leq C_1\exp\left\{-\frac{C_2\bmuj^2\eta^2}{\son \wij^2\Nij}\right\}, \quad j \in \{1,\ldots,p\}.
    \end{equation*}
\end{lma}

\begin{proof}[\textbf{\upshape Proof:}]
    Define $g(a,b,c) = ac - b^2.$  By assumptions A1 and B2, Taylor expansions of each $E[S_{jr}(t)],$ $r\in \{0,1,2\}$, yield the bound
    $$
    g\left\{E\left[\Szj(t)\right], E\left[\Soj(t)\right], E\left[\Stj(t)\right]\right\} \geq f_j^2(t)g\left\{\phi_{j0}(t),\phi_{j1}(t),\phi_{j2}(t)\right\} - 4M\bmuj(1 - \bmuj),
    $$
    where $\phi_{jr}(t) = \int_{\mathcal{U}_{tj}} u^rK(u)\du,$ $\mathcal{U}_{tj} = [-1,1]\cap \{\bmuj\inv(v - t); v \in \Tj\},$ and $M$ is defined in assumption B2.  Furthermore, suppose $U$ and $U^+$ are random variables with densities $K$ and $K^+,$ respectively, where $K^+(u) = K(u)\left\{\i01 K(u) \du\right\}^{-1}$ for $u \in [0,1].$  Then it is straightforward to show that, for any $t \in \Tj,$ $j \in \{1,\ldots,p\},$ if $\bmuj < |\Tj|/2$ for all $j$, then $g(\phi_{j0}(t),\phi_{j1}(t),\phi_{j2}(t)) \geq \tau > 0,$ where $\tau$ is the minimum of $\Var(U)$ and $\left\{\i01 K(u)\du\right\}^2\Var(U^+)$.
    By assumptions B1 and C3, there is $N$ such that $n > N$ implies $\max_{j \in \{1,\ldots,p\}} \bmuj |\Tj| < 1/2$ and $\max_{j \in \{1,\ldots,p\}} 4M\bmuj(1 - \bmuj) < m^2\tau/2,$ where $m$ is defined in assumption B2.  Thus, for $n > N,$ 
    $
    \min_{j \in \{1,\ldots,p\}} \inf_{t \in \Tj}g\left\{E\left[\Szj(t)\right], E\left[\Soj(t)\right], E\left[\Stj(t)\right]\right\} \geq m^2\tau/2.
    $
    
    Next, by continuity of $g'$, if $\max(|a|,|b|,|c|) < L$ and $\max(|a - a'|,|b-b'|,|c-c'|) < \epsilon < L,$
    then
    $
    |g(a,b,c) - g(a',b',c')| \leq 8L\epsilon.
    $
    From assumption B2,  $\normsupj{E[S_{jr}(\cdot)]} < M,$ uniformly in $n$, $j \in \{1,\ldots,p\},$ and $r \in\{0,1,2\}$.  Hence, take $\eta = m^2\tau/4$, $L = \max(\eta, M),$ and $\epsilon = \eta/8L.$  The results then follows from Lemma~\ref{lma: SjrConcIneq} since, for $n \geq N$,
    $$
    \Pr\left[\inf_{t \in \Tj} \left\{\Szj(t)\Stj(t) - \Soj^2(t)\right\} < \eta\right] \leq \Pr\left\{\max_{r \in\{ 0,1,2\}} \normsupj{S_{jr} - E[S_{jr}(\cdot)]} > \frac{\eta}{8L}\right\}.
    $$

\end{proof}

\subsubsection{\texorpdfstring{$L^2$}{L2} and Uniform Convergence of Numerator Terms}

The next two results provide exponential tail bounds related to each of the terms in \eqref{eq: MnNumTermExp1} and \eqref{eq: MnNumTermExp2} in the $L^2$ and uniform metrics, respectively.  Define $\wjns = \max_{i \in \{1,\ldots,n\}} \wij$, $\mathbb{T}_{ij} = \{\Tijl\}_{\ell = 1}^{\Nij}$, and
\begin{equation}
    \label{eq: Wijr}
    \Wij^r(t) = \wij\soNij \Kb{\bmuj}(\Tijl - t)\left(\frac{\Tijl - t}{\bmuj}\right)^r\Uijl,\quad r \in \{0,1\}. \\
\end{equation}

\begin{lma}
\label{lma: MeanNum1}
Suppose that assumptions A1, B1--B3, and C1--C3 hold.  Then there exist constants $C_1$, $C_2,$ and $C_3$ such that, for any $\epsilon > 0$, $j \in \{1,\ldots,p\}$, and $r \in \{0,1\},$
\begin{equation}
    \label{eq: MeanNum1Dis}
\begin{split}
\Pr\left(\normj{\son \Wij^r} > \epsilon\right) &\leq C_1\exp\left\{\frac{-C_2\epsilon^2}{\son \wij^2\Nij(\bmuj\inv + \Nij - 1) + \bmuj\inv \wjn\epsilon}\right\}, \\
\Pr\left(\normj{\Szj} - C_3 > \epsilon\right) &\leq C_1\exp\left\{\frac{-C_2\epsilon^2}{\bmuj\inv\son \wij^2\Nij + \bmuj^{-1/2}\wjns \epsilon}\right\}. 
\end{split}
\end{equation}
\end{lma}

\begin{proof}[\textbf{\upshape Proof:}]
    Conditional on $\mathbb{T}_{ij},$ $\Wij^r(t)$ is, for each $t \in \Tj,$ a sub-Exponential random variable with parameter at most $w_{ij}\soNij \Kb{\bmuj}(\Tijl - t)\theta^*,$ where $\theta^* = \overline{\theta} + \sigma.$  Thus, for any $\nu \geq 2,$
    $   
    E\left[\left|\Wij^r(t)\right|^\nu \mid \mathbb{T}_{ij}\right] \leq 2\nu!(2\theta^*w_{ij})^\nu\left\{\soNij \Kb{\bmuj}(\Tijl - t)\right\}^{\nu}.
    $
    Letting $K_\infty = \sup_{|u| \leq 1} K(u),$ $K_2 = \int_0^1 K^2(u)\du,$ and $M$ as in assumption B2, it follows that
    \begin{equation}
    \label{eq: Wij_nu_bound}
    E\left[\left|\Wij^r(t)\right|^\nu\right] \leq 2M\nu! (2\theta^*\wij)^\nu(K_\infty \bmuj\inv \Nij)^{\nu - 2}\left\{K_2\Nij\bmuj\inv + M\Nij(\Nij - 1)\right\}.
    \end{equation}
    By Jensen's inequality, 
    $\son E\left[\normj{\Wij^r}^\nu\right] \leq \nu L_{1nj}^2L_{2nj}^{\nu-2}/2$, where
       $L_{1nj}^2 = 16M{\theta^*}^2|\Tj|\max(K_2,M)\son \wij^2\Nij(\bmuj\inv + \Nij - 1)$ and
        $L_{2nj} = 2K_\infty\theta^*|\Tj|^{1/2}\bmuj\inv\wjn$.
    Hence, Theroem 2.5 of \cite{bosq:00} implies the first line of \eqref{eq: MeanNum1Dis}.
    
    Finally, direct calculations show that $\normsupj{E[\Szj(\cdot)]} \leq M$ and $\normsupj{\Var[\Szj(\cdot)]} \leq MK_2 \bmuj\inv \son \wij^2\Nij,$ whence
    $
    E\left[\normj{\Szj}^2\right] \leq |\Tj|\left(K_2M\bmuj\inv \son \wij^2\Nij + M^2\right) \leq C_3^2
    $
    for some fixed constant $C_3$ independent of $j,$ $p$ and $n$.  Now, the functions $\wij\Kb{\bmuj}(\Tijl - \cdot)$ are independent across both $i$ and $\ell$.  Moreover, $\normj{\wij\Kb{\bmuj}(\Tijl - \cdot)} \leq K_2^{1/2}\bmuj^{-1/2}\wjns$ for any $i$ and $\ell$, and 
    $
    \son \soNij E\left(\normj{\wij\Kb{\bmuj}(\Tijl - \cdot)}^2\right) \leq K_2M|\Tj|\bmuj\inv \son \wij^2\Nij
    $
    follow from assumptions A1 and B2.  Thus, the second line in \eqref{eq: MeanNum1Dis} follows by applying Theorem 2.6 of \cite{bosq:00}.
\end{proof}

\begin{lma}
\label{lma: MeanNum1Sup}
Suppose that assumptions A1, B1, B2, and C1--C3 hold.  Then there exist constants $C_1, C_ 2, C_3 > 0$ such that, for any $\epsilon > 0$ and $j \in \{1,\ldots,p\}$,
\begin{equation}
    \label{eq: MeanNum1SupDis}
\begin{split}
\Pr\left(\left|\son \wij \soNij \left\{|\Uijl| - E\left[|\Uijl|\right]\right\}\right| > \epsilon\right) &\leq C_1\exp\left\{\frac{-C_2\epsilon^2}{\son \wij^2\Nij^2 + \wjn\epsilon}\right\}, \\
\Pr\left(\normsupj{\Szj} - C_3 > \epsilon\right) &\leq C_1\exp\left\{ \frac{-C_2\bmuj^2\epsilon^2}{\son\wij^2\Nij}\right\}.
\end{split}
\end{equation}
\end{lma}

\begin{proof}[\textbf{\upshape Proof:}]
Conditional on $\mathbb{T}_{ij},$ $\Uijl$ are sub-Exponential random variables with parameter at most $\theta^* = \overline{\theta} + \sigma.$ Hence, for any $\nu \geq 2,$ by applying Jensen's inequality,
\begin{equation*}
    \son E\left[\left(\wij \soNij |\Uijl|\right)^\nu\right] \leq \son \wij^\nu\Nij^\nu\left(\frac{1}{\Nij}\soNij E\left[|\Uijl|^\nu\right]\right) \leq 
    \frac{\nu!}{2}\left(2\theta^*\wjn\right)^{\nu - 2}\left(16{\theta^*}^2\son \wij^2\Nij^2\right).
\end{equation*}
Then, by Theorem 2.5 of \cite{bosq:00}, the first line of \eqref{eq: MeanNum1SupDis} is established.  Finally, as previously observed, $E\left[\Szj(t)\right] \leq M,$ so one may take $C_3 = M$.  The third line of \eqref{eq: MeanNum1SupDis} then follows by applying the result of Lemma~\ref{lma: SjrConcIneq}.
\end{proof}

%% file: MeanTheoremsAndCorollariesJMVA.tex
\begin{proof}[\textbf{\upshape Proof of Theorem~\ref{thm: meanGen}:}]

The proof will be given for an FR design.  For an SR design, the proof is simplified by the fact that $S_{jr}$ does not depend on $j$, thus necessitating less stringent requirements on the bandwidth as outlined in assumption C3.  From \eqref{eq: meanDiff},
\begin{equation*}
    \normj{\mu_j - \hat{\mu}_j} \leq \left[\inf_{t \in\Tj} \left\{S_{j0}(t)S_{j2}(t) - S_{j1}^2(t)\right\}\right]\inv \left(\normsupj{S_{j2}}\normj{R_{j0} - \muj S_{j0} - \bmuj\muj'S_{j1}} + \normsupj{S_{j1}}\normj{R_{j1} - \muj S_{j1} - \bmuj \muj' S_{j2}}\right).  
\end{equation*}
Let $C_1,C_2,C_3, N,$ and $\eta$ satisfy the results of Lemmas~\ref{lma: InfDen} and \ref{lma: MeanNum1} simultaneously.  Then Lemma~\ref{lma: InfDen} implies
\begin{equation*}
\begin{split}    
&\Pr\left[\min_{j \in \{1,\ldots,p\}}\inf_{t \in \Tj}\left\{S_{j0}(t)S_{j2}(t) - S_{j1}^2(t)\right\} < \eta\right] 
    \leq p \max_{j \in \{1,\ldots,p\}} \Pr\left(\inf_{t \in \Tj}\left\{S_{j0}(t)S_{j2}(t) - S_{j1}^2(t)\right\} < \eta\right) \\
    &\hspace{1.5cm} \leq C_1p\max_{j \in \{1,\ldots,p\}}\exp\left\{-\frac{C_2\bmuj^2\eta^2}{\son \wij^2\Nij}\right\}= C_1p^{1 - C_2\eta^2\left\{\log(p)\max_{j \in \{1,\ldots,p\}}\bmuj^{-2}\son \wij^2\Nij\right\}\inv},
\end{split}
\end{equation*}
which converges to $0$ as $n \rightarrow \infty$ by assumption C3.  Similarly, since $\max_{j \in \{1,\ldots,p\}} \normsupj{E[S_{jr}(\cdot)]} < \infty$ for $r \in \{0,1,2\}$, Lemma~\ref{lma: SjrConcIneq} ensures the existence of some $R > 0$ such that
\[    
\Pr\left(\max_{j \in \{1,\ldots,p\}} \normsupj{S_{jr}} > R\right)  \leq C_1p^{1 - C_2R^2\left\{\log(p)\max_{j = 1,\ldots,p}\bmuj^{-2}\son \wij^2\Nij\right\}\inv},
\]
which again converges to $0$ as $n \rightarrow \infty$ by (C3).  Hence, 
$
\max_{j \in \{1,\ldots,p\}}\left[\inf_{t \in \Tj}\left\{S_{j0}(t)S_{j2}(t) - S_{j1}^2(t)\right\}\right]\inv = O_P(1)
$
and $\max_{j \in \{1,\ldots,p\}}\normsupj{S_{jr}} = O_P(1),$ $r \in \{0,1,2\}.$

Next, let $R > 0$ again be arbitrary and define 
\begin{equation*}
\begin{split}
        a_{n1} &= \max_{j \in \{1,\ldots,p\}} \left[\left\{\son \wij^2\Nij(\bmuj\inv + \Nij - 1)\right\}^{1/2} + \{\log(p)\}^{1/2}\bmuj\inv\wjn\right], \\ 
        a_{n2} &= \max_{j \in \{1,\ldots,p\}} \left[\left(\bmuj\inv\son\wij^2\Nij\right)^{1/2} + \{\log(p)\}^{1/2}\bmuj^{-1/2}\wjns \right].
\end{split}
\end{equation*}
Then the union bound and Lemma~\ref{lma: MeanNum1} imply that
\begin{equation}
    \label{eq: MeanNum1Rate}
\begin{split}    
        &\Pr\left(\max_{j \in \{1,\ldots,p\}} \normj{\son\Wij^0}  > Ra_{n1}\{\log(p)\}^{1/2}\right)  \leq p\max_{j \in \{1,\ldots,p\}} \Pr\left(\normj{\son \Wij^0} > Ra_{n1}\{\log(p)\}^{1/2}\right) \\
        &\hspace{1.5cm} \leq C_1p\max_{j\in \{1,\ldots,p\}}\exp\left\{-\frac{R^2a_{n1}^2\log(p)}{\son \wij^2\Nij(\bmuj\inv + \Nij - 1) + \bmuj\inv\wjn Ra_{n1}\{\log(p)\}^{1/2}}\right\} \leq C_1p^{1 - \frac{C_2^2R^2}{1 + R}}.
\end{split}
\end{equation}
By similar reasoning,
\[
\Pr\left(\max_{j \in \{1,\ldots,p\}} \normj{S_{j0}} - C_3 > Ra_{n2}\{\log(p)\}^{1/2}\right) \leq C_1p^{1 - \frac{C_2^2R^2}{1 + R}}.
\]
In addition, using assumptions A1 and B3, there is a constant $C_4$ such that
$$
\left| \son\wij\soNij K_{\bmuj}(\Tijl - t)(\Tijl - t)^2 \int_0^1 \mu_j''(t _ v(\Tijl - t))(1 - v)\mathrm{d}v \right| \leq C_4\bmuj^2 S_{j0}(t).
$$
Thus, by the definition of $\Wij^0$ and \eqref{eq: MeanNum1Rate},
\begin{equation}   
\label{eq: maxMeanNum1Rate}
\begin{split}
\max_{j \in \{1,\ldots,p\}} \normj{R_{j0} - \muj S_{j0} - \bmuj\muj'S_{j1}} &\leq \max_{j \in \{1,\ldots,p\}} \normj{\Wij^0} + \max_{j \in \{1,\ldots,p\}} \bmuj^2\normj{S_{j0}} \\
&=O_P\left(\{\log(p)\}^{1/2}a_{n1} + \left[C_3 + \{\log(p)\}^{1/2}a_{n2}\right]\max_{j \in \{1,\ldots,p\}} \bmuj^2\right).
\end{split}
\end{equation}
Since $\max_{j \in \{1,\ldots,p\}} \bmuj \rightarrow 0$ by assumption C3 and $\wjns \leq \wjn,$ it follows that
\[
\{\log(p)\}^{1/2}a_{n2} \leq \left(\log(p)\max_{j \in \{1,\ldots,p\}} \bmuj\inv \son \wij^2\Nij\right)^{1/2} + \log(p)\max_{j \in \{1,\ldots,p\}} \bmuj\inv\wjn \rightarrow 0.
\]
Thus, the rate in \eqref{eq: maxMeanNum1Rate} matches the one given in the theorem statement.  The same argument shows that \newline $\max_{j \in \{1,\ldots,p\}} \normj{R_{j1} - \muj S_{j1} - \bmuj \muj'S_{j2}}$ has the same rate, completing the proof.
\end{proof}

\begin{proof}[\textbf{\upshape Proof of Corollary~\ref{cor: meanObsSubjL2}:}]
For the OBS scheme, \mbox{$\son \wij^2\Nij(\bmuj\inv + \Nij - 1) = (n\Njbar\bmuj)\inv + (\Njtbar - \Njbar)n\inv\left(\Njbar\right)^{-2}$} and \\ $\bmuj\inv\wjn = \Njp (n\Njbar\bmuj)\inv.$
For the SUBJ scheme, $\son \wij^2\Nij(\bmuj\inv + \Nij - 1) = (n\NjH\bmuj)\inv + n\inv\left(1 - \left\{\NjH\right\}\inv\right)$ and \\ $\bmuj\inv\wjn = (n\bmuj)\inv.$ Hence, the rates given in the corollary statement follow immediately from the general rate of Theorem~\ref{thm: meanGen}.
\end{proof}

\begin{proof}[\textbf{\upshape Proof of Corollary~\ref{cor: meanOptL2}:}]
The result is given for part (a) of (i).  The remaining rates can be derived in a similar manner, and the details are omitted.  Under the conditions given in part (a), $\max_{j \in \{1,\ldots,p\}} \bmuj^2 \preceq \left\{\log(p)(n\Nmin)\inv\right\}^{2/5},$ $\min_{j \in \{1,\ldots,p\}} \bmuj\Njbar \succeq \left\{\Nmin^4\log(p)n\inv\right\}^{1/5},$ and $\min_{j \in \{1,\ldots,p\}} \bmuj \succeq \left\{\log(p)(n\Nmax)\inv\right\}^{1/5}.$
The given rate then follows by the condition on $\Nmax.$
\end{proof}

\begin{proof}[\textbf{\upshape Proof of Corollary~\ref{cor: meanOptL2_sim}:}]
The proof follows the same logic as the proof of Corollary~\ref{cor: meanObsSubjL2} and the details are omitted.
\end{proof}

\begin{proof}[\textbf{\upshape Proof of Theorem~\ref{thm: meanSupGen}:}]
Again, the proof is only given for the case of an FR design for the same reason stated at the beginning of the proof of Theorem~\ref{thm: meanGen}.  From \eqref{eq: meanDiff},
\begin{equation*}
\begin{split}
    \normsupj{\mu_j - \hat{\mu}_j} \leq \left[\inf_{t \in\Tj} \left\{S_{j0}(t)S_{j2}(t) - S_{j1}^2(t)\right\}\right]\inv &\Bigg(\normsupj{S_{j2}}\normsupj{R_{j0} - \muj S_{j0} - \bmuj\muj'S_{j1}} \\ &\hspace{1cm} + \normsupj{S_{j1}}\normsupj{R_{j1} - \muj S_{j1} - \bmuj \muj' S_{j2}}\Bigg).    
\end{split}
\end{equation*}
The same arguments used in the proof of Theorem~\ref{thm: meanGen} imply that $\max_{j \in \{1,\ldots,p\}}\left[\inf_{t \in \Tj}\left\{S_{j0}(t)S_{j2}(t) - S_{j1}^2(t)\right\}\right]\inv = O_P(1)$
and $\max_{j \in \{1,\ldots,p\}}\normsupj{S_{jr}} = O_P(1),$ $r \in \{0,1,2\}.$  Thus, the uniform rates of the terms in \eqref{eq: MnNumTermExp1} and \eqref{eq: MnNumTermExp2} are sufficient.

For any $\delta > 0$ and for each $j$, let $\chi_{j}(\delta)$ be a discrete uniform grid for $\Tj$ with spacing at most $\delta,$ and let $L$ be a constant, guaranteed by assumption B1, so that $\max_{j \in \{1,\ldots,p\}} |\chi_j(\delta)| \leq L\delta\inv,$ with $|\chi_j(\delta)|$ denoting the number of elements in the grid.  Given the decomposition in \eqref{eq: MnNumTermExp1}, consider the uniform convergence of the each term on the right-hand side of this equation.  Beginning with the first term,
\begin{equation}
    \label{eq: ProcUnif}
    \normsupj{\son \Wij^r} \leq \sup_{t \in \chi_j(\delta)}\left|\son \Wij^r(t)\right| + \sup_{|s - t| \leq \delta}\left|\son \left\{\Wij^r(s) - \Wij^r(t)\right\}\right|,
\end{equation}
where $\Wij^r$ is defined in \eqref{eq: Wijr} for $r \in \{0,1\}$.  By \eqref{eq: Wij_nu_bound}, for any $\nu \geq 2,$ $j \in \{1,\ldots,p\}$ and $t \in \Tj,$ $\son E\left[|\Wij^r(t)|^\nu\right] \leq \frac{\nu!}{2}L_{1nj}^2L_{2nj}^{\nu - 2},$ where $L_{1nj}^2 = 16M{\theta^*}^2\max(K_2, M)\son \wij^2\Nij(\bmuj\inv + \Nij - 1)$ and $L_{2nj} = 2K_\infty\theta^*\bmuj\inv\wjn.$ Thus, there exist constants $C_1$ and $C_2$ such that, for any $\epsilon,\delta > 0,$
\begin{equation}
    \label{eq: WijGridBound}
    \Pr\left\{\max_{t \in \chi_j(\delta)} \left| \son \Wij^r(t)\right| > \epsilon\right\} \leq L\delta\inv C_1\exp\left\{-\frac{C_2\epsilon^2}{\son\wij^2\Nij(\bmuj\inv + \Nij - 1) + \bmuj\inv\wjn\epsilon}\right\}.
\end{equation}
Putting $a_{n1}(\delta) = \max_{j \in \{1,\ldots,p\}} \left[\left\{\son\wij^2\Nij(\bmuj\inv + \Nij - 1)\right\}^{1/2} + \left\{\log(p\delta\inv)\right\}^{1/2}\bmuj\inv \wjn\right],$ for any $R > 0,$ \eqref{eq: WijGridBound} implies
\begin{equation}
    \label{eq: WijGridRate}
    \Pr\left[\max_{j \in \{1,\ldots,p\}}\max_{t \in \chi_j(\delta)} \left| \son \Wij^r(t)\right| > Ra_{n1}(\delta)\{\log(p\delta\inv)\}^{1/2}\right] \leq C_1L\left(p\delta\inv\right)^{1 - \frac{C_2R^2}{1 + R}}.
\end{equation}
With $\delta = \delta_n = n^{-\alpha},$ the Borel-Cantelli lemma and \eqref{eq: WijGridRate}, imply that, almost surely,
\begin{equation}
    \label{eq: WijGridBigO}
        \max_{j \in \{1,\ldots,p\}} \max_{t \in \chi_j(\delta_n)} \left| \son \Wij^r(t)\right| = 
        O\left(\max_{j \in \{1,\ldots,p\}}\left[\left\{\log(p \vee n)\son\wij^2\Nij(\bmuj\inv + \Nij - 1)\right\}^{1/2} + \log(p \vee n)\bmuj\inv \wjn\right]\right).
\end{equation}

Next, with $L_K$ being the Lipschitz constant in assumption A2 and 
$
C_4 = \limsup_{n \rightarrow \infty} \max_{j \in \{1,\ldots,p\}} E\left[|\Uijl|\right] <  \infty,
$
it follows that
\[
\sup_{|s - t| \leq \delta}\left|\son \left\{\Wij^r(s) - \Wij^r(t)\right\}\right| \leq L_K\delta\bmuj^{-2}\son \wij \soNij \left|\Uijl\right| \leq L_K\delta\bmuj^{-2}\left\{C_4 + \son \wij \soNij\left(\left|\Uijl\right| - E\left[\left|\Uijl\right|\right]\right)\right\}.
\]
Set $a_{n2} = \max_{j \in \{1,\ldots,p\}}\left[\left(\son \wij^2\Nij^2\right)^{1/2} + \left\{\log(p)\right\}^{1/2}\wjn\right]$.
Then Lemma~\ref{lma: MeanNum1Sup} implies that there are constants $C_1, C_2 > 0$ such that, for any $R > 0,$
\begin{equation*}
    \Pr\left[\max_{j \in \{1,\ldots,p\}}\left|\son \wij \soNij\left(|\Uijl| - E\left[\left|\Uijl\right|\right]\right)\right| > Ra_{n2}\left\{\log(p)\right\}^{1/2}\right] \leq C_1p^{1 - \frac{C_2R^2}{1 + R}}.
\end{equation*}
By the Borel-Cantelli Lemma, assumption C4, and again taking $\delta = \delta_n = n^{-\alpha},$ it follows that, almost surely,
\begin{equation}
    \label{eq: WijDiffBigO}
    \begin{split}
    \max_{j \in \{1,\ldots,p\}} \sup_{|s - t| < \delta_n}\left| \son \{\Wij^r(s) - \Wij^r(t)\}\right| &= O\left(n^{-\alpha}\max_{j \in \{1,\ldots,p\}} \bmuj^{-2}\right) \\
    &= o\left(\max_{j \in \{1,\ldots,p\}}\left[\left\{\log(p \vee n)\son\wij^2\Nij(\bmuj\inv + \Nij - 1)\right\}^{1/2} + \log(p\vee n)\bmuj\inv\wjn\right]\right).
    \end{split}
\end{equation}
Then \eqref{eq: ProcUnif}, \eqref{eq: WijGridBigO} and \eqref{eq: WijDiffBigO} together imply the result since, almost surely,
\begin{equation*}
\max_{j \in \{1,\ldots,p\}} \normsupj{\son \Wij^r} = O\left(\max_{j \in \{1,\ldots,p\}}\left[\left\{\log(p \vee n)\son\wij^2\Nij(\bmuj\inv + \Nij - 1)\right\}^{1/2} + \log(p\vee n)\bmuj\inv\wjn\right]\right).
\end{equation*}
 
\end{proof}

\begin{proof}[\textbf{\upshape Proofs of Corollaries~\ref{cor: meanObsSubjSup}--\ref{cor: meanOptSupSim}:}]
The proofs are similar to those of Corollaries~\ref{cor: meanObsSubjL2}--\ref{cor: meanOptL2_sim}, and are omitted.
\end{proof}

%% file: CovLemmasJMVA.tex
For $q,r\in\{0,1,2\},$ and $j,k\in\{1,\dots, p\},$ define
\begin{equation}
\label{eq: SjkqrRjkqrDefs}
\begin{split}
    S_{jkqr}(s,t) &= \son \vijk \sIijk \Kb{\bgammaj}\pr{\Tijl-s}\Kb{\bgammak}\pr{\Tikm - t}\pr{\frac{\Tijl-s}{\bgammaj}}^q\pr{\frac{\Tikm - t}{\bgammak}}^r, \\
    R_{jkqr}(s,t) &= \son \vijk \sIijk \Kb{\bgammaj}\pr{\Tijl-s}\Kb{\bgammak}\pr{\Tikm - t}\pr{\frac{\Tijl-s}{\bgammaj}}^q\pr{\frac{\Tikm - t}{\bgammak}}^r\Zijklm.
\end{split}
\end{equation}
Additionally, dropping the functional arguments $s$ and $t$ for the component elements defined above in \eqref{eq: SjkqrRjkqrDefs}, set
\begin{equation}
    \label{eq: Qjkr}
    Q_{jk0} = \Sjktz\Sjkzt - \Sjkoo^2,\quad
    Q_{jk1} = \Sjkoz\Sjkzt - \Sjkzo\Sjkoo, \quad
    Q_{jk2} = \Sjkoz\Sjkoo - \Sjkzo\Sjktz.
\end{equation}
As mentioned in Section~\ref{sec: cov}, the mean is assumed to be known; that is, $\hgammajk(s,t) = \hbz,$ where $\hbz$ is computed as in \eqref{eq: CovEst} with $\Zijklm = \{\Yijl - \muj(\Tijl)\}\{\Yikm - \muk(\Tikm)\}$ instead of the true empirical version in which the $\Yijl$ and $\Yikm$ are centered with respect to their estimated means.  Then the error in covariance estimation can be expressed as
\begin{equation}
    \label{eq: covDiff}
    \begin{split}
    \hgammajk(s,t) - \gammajk(s,t) &= \left(Q_{jk0}\Sjkzz - Q_{jk1}\Sjkoz + Q_{jk2}\Sjkzo\right)\inv \\
    &\hspace{0.4cm}\times \left[Q_{jk0}\left\{\Rjkzz - \gammajk(s,t)\Sjkzz - \bgammaj\pgjks(s,t)\Sjkoz - \bgammak\pgjkt(s,t)\Sjkzo\right\} \right.\\
    &\hspace{0.6cm} - Q_{jk1}\left\{\Rjkoz - \gammajk(s,t)\Sjkoz - \bgammaj\pgjks(s,t)\Sjktz - \bgammak\pgjkt(s,t)\Sjkoo\right\} \\
    &\hspace{0.6cm} + \left.Q_{jk2}\left\{\Rjkzo - \gammajk(s,t)\Sjkzo - \bgammaj\pgjks(s,t)\Sjkoo - \bgammak\pgjkt(s,t)\Sjkzt\right\}\right]. \\
    \end{split}
\end{equation}

The terms in square brackets in the last three lines of \eqref{eq: covDiff} can be broken down further as follows.  Define
\begin{equation*}
    \Uijklm = \Uij(\Tijl)\Uik(\Tikm) - \gammajk(s,t) + \Uij(\Tijl)\epsikm + \Uik(\Tikm)\epsijl + \epsijl\epsikm
\end{equation*}
and define $\tilde{K}_b(u) = \Kb{b}(u)(u/b).$  Furthermore, set
\begin{equation*}
    B_{jk\ell m}(s,t) = \gammajk(\Tijl,\Tikm) - \gammajk(s,t) - \bgammaj(\Tijl - s)\pgjks(s,t) - \bgammak(\Tijl - t)\pgjkt(s,t).
\end{equation*}
Then
\begin{equation}
    \label{eq: CovNumTermExp1}
    \begin{split}
    \Rjkzz - \gammajk(s,t)\Sjkzz - \bgammaj\pgjks\Sjkoz - \bgammak\pgjkt\Sjkzo &= \son \vijk \sIijk \Kb{\bgammaj}(\Tijl - s)\Kb{\bgammak}(\Tikm - t)\Uijklm \\
    & \hspace{1cm} + \son \vijk \sIijk \Kb{\bgammaj}(\Tijl - s)\Kb{\bgammak}(\Tikm - t)B_{jk\ell m}(s,t),
    \end{split}    
\end{equation}
\begin{equation}
    \label{eq: CovNumTermExp2}
    \begin{split}
    \Rjkoz - \gammajk(s,t)\Sjkoz - \bgammaj\pgjks\Sjktz - \bgammak\pgjkt\Sjkoo &= \son \vijk \sIijk \tKb{\bgammaj}(\Tijl - s)\Kb{\bgammak}(\Tikm - t)\Uijklm \\
    & \hspace{1cm} + \son \vijk \sIijk \tKb{\bgammaj}(\Tijl - s)\Kb{\bgammak}(\Tikm - t)B_{jk\ell m}(s,t),
    \end{split} 
\end{equation}
\begin{equation}
    \label{eq: CovNumTermExp3}
    \begin{split}
    \Rjkzo - \gammajk(s,t)\Sjkzo - \bgammaj\pgjks\Sjkoo - \bgammak\pgjkt\Sjkzt &=\son \vijk \sIijk \Kb{\bgammaj}(\Tijl - s)\tKb{\bgammak}(\Tikm - t)\Uijklm \\
    & \hspace{1cm} + \son \vijk \sIijk \Kb{\bgammaj}(\Tijl - s)\tKb{\bgammak}(\Tikm - t)B_{jk\ell m}(s,t).
    \end{split} 
\end{equation}

\subsubsection{Uniform Convergence of \texorpdfstring{$S_{jkqr}$}{Sjkqr}, \texorpdfstring{$q, r \in \{0, 1, 2\}$ and $0 \leq q + r \leq 2$}{0 <= q + r <= 2}}
\label{ss: Sjkqr}

To begin, a result similar to Lemma~\ref{lma: SjrConcIneq} will be established, however the proof is much more involved.  Specifically, an analog of Lemma 1.1 of \cite{marc:84} is proved that is suitable for the quantities $S_{jkqr}$ in \eqref{eq: SjkqrRjkqrDefs} that are not simply weighted kernel density estimators.  

\begin{lma}
\label{lma: MarcZinnXCov}
Let $G_1$ and $G_2$ be arbitrary cumulative distribution functions.  For any $n \in \mathbb{N},$ let $N_{ij} \in \mathbb{N},$ $i \in \{1,\ldots, n\}$, $j \in \{1,2\},$ be arbitrary and consider independent arrays of random variables $\mathbb{T}_j = \left\{\Tijl:\, \ell \in \{1,\ldots,\Nij\},\, i \in \{1,\ldots,n\}\right\},$ where the elements of $\mathbb{T}_j$ are independent and identically distributed according to $G_j,$ $j \in \{1,2\}.$  Furthermore, let $v_i$, $i \in \{1, \ldots, n\},$ be arbitrary constants and define $\omega_1^2 = \son v_i^2N_{i1}^2N_{i2},$ $\omega_2^2 = \son v_i^2N_{i1}N_{i2}^2$, and $\overline{\omega} = \max_{j \in \{1,2\}} \omega_j.$  Set
$
\mathbb{A}_{n}(s,t) = \son v_i\sum_{\ell = 1}^{N_{i1}}\sum_{m = 1}^{N_{i2}}\left\{ \mathbf{1}(T_{i1\ell} \leq s) \mathbf{1}(T_{i2m} \leq t) - G_1(s)G_2(t)\right\}.
$
Then, for any $\epsilon > 0,$
\[
\Pr\left\{\sup_{s,t \in \mathbb{R}} |\mathbb{A}_n(s,t)| > \epsilon\overline{\omega}\right\} \leq \left(3 + \frac{3\sqrt{2\pi}\epsilon}{2} + \frac{2\pi\epsilon^2}{9}\right)e^{-\epsilon^2/288}.
\]
\end{lma}

\begin{proof}[\textbf{\upshape Proof:}]
First, note that
$    \mathbb{A}_n(s,t) = \mathbb{A}_{n0}(s,t) + \mathbb{A}_{n1}(s,t) + \mathbb{A}_{n2}(s,t),$
where
\begin{equation*}
    \label{eq:AnComponents}
    \begin{split}
        \mathbb{A}_{n0}(s,t) &= \son v_i \sum_{\ell = 1}^{N_{i1}}\sum_{m = 1}^{N_{i2}} \left\{\mathbf{1}(T_{i1\ell} \leq s) - G_1(s)\right\}\left\{\mathbf{1}(T_{i2m} \leq s) - G_2(t)\right\}, \\
        \mathbb{A}_{n1}(s,t) &= G_2(t)\son v_iN_{i2}\sum_{\ell = 1}^{N_{i1}} \left\{\mathbf{1}(T_{i1\ell} \leq s) - G_1(s)\right\}, \\        
        \mathbb{A}_{n2}(s,t) &= G_1(s)\son v_iN_{i1}\sum_{m = 1}^{N_{i2}} \left\{\mathbf{1}(T_{i2m} \leq t) - G_2(t)\right\}.
    \end{split}
\end{equation*}
As $\sup_{s}|G_j(s)| = 1$ for $j \in \{1,2\},$ it follows immediately from Lemma 1.1 of \cite{marc:84} that, for any $\epsilon > 0,$
\begin{equation}
    \label{eq:An12ProbBound}
    \Pr\left\{\sup_{s,t} \left|\mathbb{A}_{nj}(s,t)\right| > \epsilon\omega_j\right\} \leq \left(1 + 2\sqrt{2\pi}\epsilon\right)e^{-\epsilon^2/8}, \quad j \in \{1,2\}.
\end{equation}
To uniformly bound $\mathbb{A}_{n0},$ it will be established that, for any $\lambda > 0$ and $\underline{\omega} = \min_{j \in \{1,2\}} \omega_j,$
\begin{equation}
    \label{eq:An0MGFBound}
    E\left[\exp\left\{\lambda \sup_{s,t} \left|\mathbb{A}_{n0}(s,t)\right|\right\}\right] \leq 1 + 64\sqrt{2\pi}\lambda \underline{\omega} + 16\sqrt{2\pi}\lambda \underline{\omega}\left(1 + 32\sqrt{2\pi}\lambda \underline{\omega}\right)e^{32\lambda^2\underline{\omega}^2}.
\end{equation}
Once established, \eqref{eq:An0MGFBound} implies that, for any $\epsilon > 0$ and $\lambda = \epsilon(64\underline{\omega})\inv,$
\begin{equation}
    \label{eq:An0ProbBound}
        \Pr\left\{\sup_{s,t} \left|\mathbb{A}_{n0}(s,t)\right| > \epsilon \underline{\omega}\right\} \leq e^{-\lambda\underline{\omega}\epsilon} E\left[\exp\left\{\lambda \sup_{s,t} \mathbb{A}_{n0}(s,t)\right\}\right] 
        \leq \left(1 + \frac{5\sqrt{2\pi}\epsilon}{4} + \frac{\pi\epsilon^2}{2}\right)e^{-\epsilon^2/128}.
\end{equation}
Together, \eqref{eq:An12ProbBound} and \eqref{eq:An0ProbBound} imply that
\begin{equation*}
\begin{split}
    &\Pr\left\{\sup_{s,t \in \mathbb{R}} |\mathbb{A}_n(s,t)| > \epsilon\overline{\omega}\right\} \leq \Pr\left\{\sup_{s,t} \left|\mathbb{A}_{n0}(s,t)\right| > \frac{2\epsilon \underline{\omega}}{3}\right\} + \sum_{j = 1}^2 \Pr\left\{\sup_{s,t} \left|\mathbb{A}_{nj}(s,t)\right| > \frac{\epsilon \omega_j}{6}\right\} \\
    &\hspace{1.5cm} \leq \left(1 + \frac{5\sqrt{2\pi}\epsilon}{6} + \frac{2\pi\epsilon^2}{9}\right)e^{-\epsilon^2/288} + 2\left(1 + \frac{\sqrt{2\pi}\epsilon}{3}\right)e^{-\epsilon^2/288}
    = \left(3 + \frac{3\sqrt{2\pi}\epsilon}{2} + \frac{2\pi\epsilon^2}{9}\right)e^{-\epsilon^2/288}.
\end{split}
\end{equation*}
Thus, it remains only to prove \eqref{eq:An0MGFBound}.

Let $(\mathbb{T}_1', \mathbb{T}_2')$ represent an independent copy of $(\mathbb{T}_1,\mathbb{T}_2)$, leading to $\mathbb{A}_{n0}'(s,t)$ in analogy to $\mathbb{A}_{n0}(s,t).$  Let $e_{i1\ell}$ and $e_{i2m},$ $i \in \{1, \ldots, n\},$ $\ell \in \{1,\ldots,N_{i1}\},$ $m \in \{1,\ldots,N_{i2}\}$ denote independent and identically distributed Radamacher variables, that is, $\Pr(e_{i1\ell} = 1) = \Pr(e_{i1\ell} = -1) = 1/2$.  By symmetry, simple calculations show that $\mathbb{A}_{n0}(s,t) - \mathbb{A}_{n0}'(s,t) \overset{D}{=} Y_1(s,t) + Y_2(s,t)$,
where equality of distribution holds at the process level and
\[
\begin{split}
Y_1(s,t) &= \son v_i\sum_{\ell = 1}^{N_{i1}}\sum_{m = 1}^{N_{i2}}\left\{\mathbf{1}(T_{i1\ell} \leq s) - \mathbf{1}(T_{i1\ell}' \leq s)\right\}\left\{\mathbf{1}(T_{i2m} \leq t) - G_2(t)\right\}e_{i1\ell} \\
Y_2(s,t) &= \son v_i\sum_{\ell = 1}^{N_{i1}}\sum_{m = 1}^{N_{i2}}\left\{\mathbf{1}(T_{i1\ell} \leq s) - G_1(s)\right\}\left\{\mathbf{1}(T_{i2m} \leq t) - \mathbf{1}(T_{i2m}' \leq t)\right\}e_{i2m}.
\end{split}
\]
Then, for $\lambda > 0,$ by applying Jensen's inequality, the triangle inequality, and Cauchy-Schwarz,
\begin{equation}
    \label{eq:An0MGFpt1}
    \begin{split}
        E\left[\exp\left\{\lambda \sup_{s,t} |\mathbb{A}_{n0}(s,t)|\right\}\right] &\leq E\left(\exp\left\{\sup_{s,t}|\mathbb{A}_{n0}(s,t) - \mathbb{A}_{n0}'(s,t)|\right\}\right)
        \leq E\left(\exp\{\lambda \sup_{s,t}|Y_1(s,t)|\}\exp\{\lambda \sup_{s,t}|Y_2(s,t)|\}\right) \\
        &\leq \left(E\left[\exp\left\{2\lambda \sup_{s,t}|Y_1(s,t)|\right\}\right]E\left[\exp\left\{2\lambda \sup_{s,t}|Y_2(s,t)|\right\}\right]\right)^{1/2}.
    \end{split}
\end{equation}

Next, define $v_{i\ell}(t) = v_i\sum_{m = 1}^{N_{i2}}\left\{\mathbf{1}(T_{i2m} \leq t) - G_2(t)\right\}$, let $N_1^S = \son N_{i1}$ and suppose that $\{T_{1\ell'}\}_{\ell' = 1}^{N_1^S}$ is a non-decreasing ordering of $\mathbb{T}_1.$ If $\ell',$ $i$ and $\ell$ are such that $T_{1\ell'} = T_{i1\ell},$ set $v_{\ell'}(t) = v_{i\ell}(t)$ and $e_{1\ell'} = e_{i1\ell}.$  Then applications of Cauchy-Schwarz and the triangle inequality, along with the fact that $\mathbb{T}_1 \overset{D}= \mathbb{T}_1'$, yield
\begin{equation}
    \label{eq:Y1MGFbndinit}
    \begin{split}
    E\left[\exp\left\{2\lambda \sup_{s,t} |Y_1(s,t)|\right\}\right] &\leq E\left[\exp\left\{4\lambda \sup_{s,t} \left|\son \sum_{\ell = 1}^{N_{i1}} v_{i\ell}(t)\mathbf{1}(T_{i1\ell} \leq s)e_{i1\ell}\right|\right\}\right]
    = E\left[\exp\left\{4 \lambda \max_{\ell'' \in \{1,\ldots,N_1^S\}} \sup_t \left| \sum_{\ell' = 1}^{\ell''} v_{\ell'}(t)e_{1\ell'}\right| \right\}\right] \\
    &= 1 + \int_1^\infty \Pr\left(\max_{\ell'' \in \{1,\ldots,N_1^S\}} \sup_t \left| \sum_{\ell' = 1}^{\ell''} v_{\ell'}(t)e_{1\ell'}\right| > \frac{\log(\gamma)}{4\lambda}\right) \mathrm{d}\gamma.
    \end{split}
\end{equation}
Let $\eta > 0.$  Then, by Levy's inequality, applied conditionally on $\mathbb{T}_2$, 
\begin{equation}
    \label{eq:Y1ProbBound}
    \Pr\left\{\max_{\ell'' \in \{1,\ldots,N_1^S\}} \sup_t \left| \sum_{\ell' = 1}^{\ell''} v_{\ell'}(t)e_{1\ell'}\right| > \eta\right\} \leq 2\Pr\left[\sup_t \left|\son v_i \sum_{\ell = 1}^{N_{i1}}\sum_{m = 1}^{N_{i2}}e_{i1\ell}\left\{\mathbf{1}(T_{i2m} \leq t) - G_2(t)\right\}\right| > \eta\right].
\end{equation}
To bound this probability, the moment generating function will again be bounded by symmetrization.  Specifically, for any $\delta > 0,$ apply Jensen's (conditionally on the $e_{i1\ell}$) to obtain
\begin{equation}
    \label{eq:CnMGFBound}
E\left[\exp\left\{\delta\sup_t \left|\son v_i \sum_{\ell = 1}^{N_{i1}}\sum_{m = 1}^{N_{i2}}e_{i1\ell}\left\{\mathbf{1}(T_{i2m} \leq t) - G_2(t)\right\}\right|\right\}\right]  \leq E\left[\exp\left\{ 2\delta \sup_t \left|\son v_i\sum_{\ell = 1}^{N_{i1}}\sum_{m = 1}^{N_{i2}}\mathbf{1}(T_{i2m} \leq t)e_{i1\ell}e_{i2m}\right|\right\}\right].
\end{equation}
Let $N_2^S = \sum_{i = 1}^n N_{i2}$ and let $\{T_{2m'}\}_{m' = 1}^{N_{2}^S}$ be a non-decreasing ordering of $\mathbb{T}_2.$ Define $v_{im} = v_i\sum_{\ell = 1}^{N_{i1}}e_{i1\ell}$ and, if $T_{2m'} = T_{i2m},$ set $v_{m'}' = v_{im}$ and $e_{2m'} = e_{i2m}.$ Then, for any $\tau > 0,$ another application of Levy's inequality combined with Hoeffding's inequality yields
\begin{equation}
    \label{eq:CnStarProbBound}
    \begin{split}
        &\Pr\left\{\sup_t \left|\son v_i\sum_{\ell = 1}^{N_{i1}}\sum_{m = 1}^{N_{i2}}\mathbf{1}(T_{i2m} \leq t)e_{i1\ell}e_{i2m}\right| > \tau\right\} = \Pr\left(\max_{m'' \in \{1,\ldots,N_2^S\}} \left|\sum_{m' = 1}^{m''} v_{m'}'e_{2m'}\right| > \tau\right) \\
        &\hspace{1.5cm} \leq 2\Pr\left(\left|\son v_i\sum_{\ell = 1}^{N_{i1}}\sum_{m = 1}^{N_{i2}} e_{i1\ell}e_{i2m}\right| > \tau\right) \leq 4\exp\left\{-\frac{\tau^2}{2\underline{\omega}}\right\}.
    \end{split}
\end{equation}
Combining \eqref{eq:CnMGFBound} and \eqref{eq:CnStarProbBound}, it can be concluded that 
\begin{equation*}
    \label{eq:CnMGFBound2}
    \begin{split}
        &E\left[\exp\left\{\delta\sup_t \left|\son v_i \sum_{\ell = 1}^{N_{i1}}\sum_{m = 1}^{N_{i2}}e_{i1\ell}\left\{\mathbf{1}(T_{i2m} \leq t) - G_2(t)\right\}\right|\right\}\right] 
        = 1 + \int_{1}^\infty \Pr\left(\sup_t \left| \son v_i \sum_{\ell = 1}^{N_{i1}}\sum_{m = 1}^{N_{i2}}\mathbf{1}(T_{i2m} \leq t)e_{i1\ell}e_{i2m}\right| > \frac{\log(\gamma)}{2\delta}\right) \mathrm{d}\gamma \\
        &\hspace{1.5cm} \leq 1 + 8\delta\underline{\omega}\int_0^\infty e^{2\delta\underline{\omega}\tau}e^{-\tau^2/2}\mathrm{d}\tau 
        \leq 1 + 8\delta\underline{\omega}\sqrt{2\pi} e^{2\delta^2\underline{\omega}^2},
    \end{split}
\end{equation*}
where the change of variables $\log(\gamma) = 2\delta\underline{\omega}\tau$ and the integral bound $\int_0^\infty e^{rx}e^{-x^2/2}\mathrm{d}x \leq \sqrt{2\pi}e^{r^2/2}$ for $r > 0$ have been used.  Hence, for any $\eta > 0,$ setting  $\delta = \eta/4\underline{\omega}$ yields the bound
\begin{equation}
    \label{eq:CnProbBound2}
    \begin{split}
    \Pr\left[\sup_t \left|\son v_i \sum_{\ell = 1}^{N_{i1}}\sum_{m = 1}^{N_{i2}}e_{i1\ell}\left\{\mathbf{1}(T_{i2m} \leq t) - G_2(t)\right\}\right| > \eta\underline{\omega}\right]
    &\leq e^{-\delta\underline{\omega}\eta}E\left(\exp\left\{\delta\underline{\omega}\sup_t \left|\son v_i \sum_{\ell = 1}^{N_{i1}}\sum_{m = 1}^{N_{i2}}e_{i1\ell}\left[\mathbf{1}(T_{i2m} \leq t) - G_2(t)\right]\right|\right\}\right)\\
    &\leq e^{-\eta^2/4}\left(1 + 2\eta\sqrt{2\pi}e^{\eta^2/8}\right) \leq (1 + 2\eta\sqrt{2\pi})e^{-\eta^2/8}.
    \end{split}
\end{equation}

Returning to a bound for the first expectation in \eqref{eq:Y1MGFbndinit}, \eqref{eq:Y1ProbBound} and \eqref{eq:CnProbBound2} imply
\begin{equation*}
    \label{eq:Y1MGFBound}
    \begin{split}
        &E\left[\exp\left\{2\lambda\sup_{s,t}|Y_1(s,t)|\right\}\right] \leq 1 + 2\int_1^\infty \Pr\left(\sup_t \left|\son v_i \sum_{\ell = 1}^{N_{i1}}\sum_{m = 1}^{N_{i2}}e_{i1\ell}\left[\mathbf{1}(T_{i2m} \leq t) - G_2(t)\right]\right|  > \frac{\log(\gamma)}{4\lambda}\right)\mathrm{d}\gamma \\
        &\hspace{1.5cm} \leq 1 + 16\lambda\underline{\omega}\int_1^\infty (1 + 4\sqrt{2\pi}\eta)e^{8\lambda\underline{\omega}\eta}e^{-\eta^2/2}\mathrm{d}\eta 
        \leq 1 + 16\lambda\underline{\omega}\left[4\sqrt{2\pi} + \sqrt{2\pi}\left(1 + 32\lambda\underline{\omega}\sqrt{2\pi}\right)e^{32\lambda^2\underline{\omega}^2}\right] \\
        &\hspace{1.5cm} \leq 1 + 64\sqrt{2\pi}\lambda\underline{\omega} + 16\sqrt{2\pi}\lambda\underline{\omega}\left(1 + 32\sqrt{2\pi}\lambda \underline{\omega}\right)e^{32\lambda^2\underline{\omega}^2},
    \end{split}
\end{equation*}
where the change of variable $\log(\gamma) = 8\lambda\underline{\omega}\tau$ and the integral bound $\int_0^\infty xe^{rx}e^{-x^2/2}\mathrm{d}x \leq 1 + r\sqrt{2\pi}e^{r^2/2}$ for $r > 0$ have been used.  By symmetry, the same bound applies to the term involving $Y_2$ in \eqref{eq:An0MGFpt1}, so \eqref{eq:An0MGFBound} follows.
\end{proof}

\begin{lma}
    \label{lma: MarcZinnACov}
Let $G$ be an arbitrary cumulative distribution function.  For any $n \in \mathbb{N},$ let $N_{i} \in \mathbb{N},$ $i \in \{1,\ldots, n\}$ be arbitrary and consider an independent array of random variables $\mathbb{T} = \left\{T_{i\ell}:\, \ell \in \{1,\ldots,N_i\},\, i \in \{1,\ldots,n\}\right\},$ whose elements are independent and identically distributed according to $G$.  Furthermore, let $v_i$, $i \in \{1, \ldots, n\},$ be arbitrary constants and define $\omega^2 = \son v_i^2N_i(N_i - 1)^2.$  Set
$
\mathbb{B}_{n}(s,t) = \son v_i \sum_{\ell = 1}^{N_{i}}\sum_{m \neq \ell}\left\{\mathbf{1}(T_{i\ell} \leq s)\mathbf{1}(T_{im} \leq t) - G(s)G(t)\right\}
$
Then, for any $\epsilon > 0,$
\[
\Pr\left(\sup_{s,t \in \mathbb{R}} |\mathbb{B}_n(s,t)| > \epsilon\omega \right) \leq \left(1 + 2\sqrt{2\pi}\epsilon\right)e^{-\epsilon^2/8}
\]

\end{lma}

\begin{proof}[\textbf{\upshape Proof:}]
    The proof is similar to that of Lemma~\ref{lma: MarcZinnXCov} and the details are omitted.
\end{proof}

\begin{lma}
\label{lma: SjkqrConcIneq}
Suppose that assumptions A1 and B1 hold, and recall the definition of $\omega_{jk}$ in \eqref{eq: omegajk}.  Then, there exist $C_1, C_2 > 0$ such that, for any $\epsilon > 0$, integers $q,r\in\{0,1,2\}$ such that $0 \leq q + r \leq 2$ and $j,k \in \{1,\ldots,p\},$
\[
\Pr\left\{\normsupjk{S_{jkqr}(\cdot,*) - E[S_{jkqr}(\cdot,*)]} > \epsilon\right\} \leq C_1\exp\left\{\frac{-C_2\bgammaj^2\bgammak^2\epsilon^2}{\omega_{jk}^2}\right\} 
\]
under an FR design and, under an SR design,
\[
\Pr\left\{\normsupjk{S_{jkqr}(\cdot,*) - E[S_{jkqr}(\cdot,*)]} > \epsilon\right\} \leq C_1\exp\left\{\frac{-C_2\bgamma^4\epsilon^2}{\omega^2}\right\}.
\]
\end{lma}

\begin{proof}[\textbf{\upshape Proof:}]
   The proof under both designs will be shown for $j \neq k$; the case $j = k$ follows similar, but simpler arguments.  Begin with the FR design.  Standard derivations arising from application of Riemann-Stiltjes integration by parts yield
    $
    \normsupjk{S_{jkqr}(\cdot,*) - E\left[S_{jkqr}(\cdot, *)\right]} \leq V_K\bgammaj\inv\bgammak\inv\normsupjk{\mathbb{A}_{jk}},
    $
    where $V_K$ is a constant depending only on the kernel $K$,
    $
    \mathbb{A}_{jk}(s,t) = \son \vijk \sum_{\ell = 1}^{\Nij}\sum_{m = 1}^{\Nik}\left\{\mathbf{1}(\Tijl \leq s)\mathbf{1}(T_{ijk} \leq t) - F_j(s)F_k(t)\right\},
    $
    and $F_j$ is the cdf corresponding to $f_j.$  Thus, the first inequality of the lemma follows from an application of Lemma~\ref{lma: MarcZinnXCov}.

    In the case of an SR design, write $S_{qr}$ and $v_i$ for the common values of $S_{jkqr}$ and $\vijk$.  Then \newline
    $
    \normsup{S_{qr}(\cdot,*) - E(S_{qr}(\cdot, *))} \leq V_K\bgamma^{-2}\normsup{\mathbb{B}},
    $
    where
    $
            \mathbb{B}(s,t) = \son v_i \sum_{\ell = 1}^{N_i} \sum_{m \neq l}\left[\mathbf{1}(T_{i\ell} - s)\mathbf{1}(T_{im} - t) - F(s)F(t)\right].
    $
    The result then follows by applying Lemma~\ref{lma: MarcZinnACov}.
\end{proof}

\begin{lma}
    \label{lma: InfDenCov}
    Suppose assumptions A1, B1, B2, and D3 hold.  Define 
    \begin{equation*}
        S_{jk}^*(s,t) = Q_{jk0}\Sjkzz - Q_{jk1}\Sjkoz + Q_{jk2}\Sjkzo,
    \end{equation*} 
    and define $\omega_{jk}$ and $\omega$ as in \eqref{eq: omegajk}.  Then there exist $C_1,C_2,\eta,N > 0$ such that, for any $n \geq N$ and $j,k \in \{1,\ldots,p\},$
    $$
    \Pr\left\{\inf_{(s,t) \in \Tj \times \Tk} S_{jk}^*(s,t) < \eta\right\} \leq C_1\exp\left\{\frac{-C_2\bgammaj^2\bgammak^2\epsilon^2}{\omega_{jk}^2}\right\}
    $$
    under an FR design and, under an SR design,
    $$
    \Pr\left\{\inf_{(s,t) \in \T^2} S_{jk}^*(s,t) < \eta\right\} \leq C_1\exp\left\{\frac{-C_2\bgamma^4\epsilon^2}{\omega^2}\right\}.
    $$
\end{lma}
\begin{proof}[\textbf{\upshape Proof:}]
    The proof follows the same lines as the proof of Lemma~\ref{lma: InfDen}, so only a sketch will be provided.  Define
    \[
    \begin{split}
    \overline{S}_{jk}^*(s,t) &= \left(E\left[\Sjktz\right]E\left[\Sjkzt\right] - E\left[\Sjkoo\right]^2\right)E\left[\Sjkzz\right]  - \left(E\left[\Sjkoz\right]E\left[\Sjkzt\right] - E\left[\Sjkzo\right]E\left[\Sjkoo\right]\right)E\left[\Sjkoz\right] \\
    &\hspace{1cm} + \left(E\left[\Sjkoz\right]E\left[\Sjkoo\right] - E\left[\Sjkzo\right]E\left[\Sjktz\right]\right)E\left[\Sjkzo\right].
    \end{split}
    \]
    First, one establishes that
    $
    \overline{S}_{jk}^*(s,t) \geq \tau f_j^3(s)f_k^3(t) - o(1)
    $
    uniformly in $j,k,s,$ and $t,$ where $\tau$ depends only on $K$.  Then, by assumptions B2 and D3, there exists $N$ such that, when $n > N,$
    $
    \min_{j,k \in \{1,\ldots,p\}} \inf_{(s,t) \in \Tj \times \Tk}\overline{S}_{jk}^*(s,t) \geq \frac{m^4\tau}{2}.
    $
    The remainder of the argument follows from uniform continuity and Lemma~\ref{lma: SjkqrConcIneq}.
\end{proof}

\begin{lma}
    \label{lma: SupNumCov}
    Suppose assumptions A1, B1, B2, and D3 hold.  Let $Q_{jkr}$ be as defined in \eqref{eq: Qjkr}, $r \in \{0,1,2\},$ and let $\omega_{jk}$ be as in \eqref{eq: omegajk}.  Then there exist $C_1,C_2,C_3 > 0$ such that, for any $\epsilon > 0$, $j,k \in \{1,\ldots,p\},$ and $r \in \{0,1,2\},$
    $$
    \Pr\left(\normsupjk{Q_{jkr}} - C_3 > \epsilon \right) \leq C_1\exp\left\{\frac{-C_2\bgammaj^2\bgammak^2\epsilon^2}{\omega_{jk}^2}\right\}
    $$
    under an FR design and, under an SR design,
    $$
    \Pr\left(\normsupjk{Q_{jkr}} - C_3 > \epsilon\right) \leq C_1\exp\left\{\frac{-C_2\bgamma^4\epsilon^2}{\omega^2}\right\}.
    $$
\end{lma}
\begin{proof}[\textbf{\upshape Proof:}]
    Assumptions B2 and D3 imply that there is $C_3$ such that $\normsupjk{E(Q_{jkr})} \leq C_3.$ Apply Lemma~\ref{lma: SjkqrConcIneq}.
\end{proof}

\subsubsection{\texorpdfstring{$L^2$}{L2} and Uniform Convergence of Numerator Terms}

To derive the $L^2$ convergence rates in Section~\ref{sec: cov}, preliminary concentration inequalities will be provided for the $L^2$ norms of the quantities in \eqref{eq: CovNumTermExp1}--\eqref{eq: CovNumTermExp3}.  For $(q,r) \in \{(0, 0), (1,0), (0,1)\},$ define $\mathbb{T}_{ijk} = \left\{(\Tijl, \Tikm): (\ell, m) \in \Iijk\right\}$ and
\begin{equation}
    \label{eq: Wijkqr}
    \Wijk^{qr}(s,t) = \vijk\sIijk \Kb{\bgammaj}(\Tijl - s)\Kb{\bgammak}(\Tikm - t)\left(\frac{\Tijl - s}{\bgammaj}\right)^q\left(\frac{\Tikm - t}{\bgammak}\right)^r \Uijklm
\end{equation}
Recall the definitions of the rates $\qjkn$ in \eqref{eq: qjkn} and define $\vjkns = \max_{i \in \{1,\ldots,n\}} \vijk.$  As $|B_{jk}(s, t)| \leq C\bgammaj^2\bgammak^2$ for a universal constant $C$ by assumption B4, the term on the last line of each of \eqref{eq: CovNumTermExp1}--\eqref{eq: CovNumTermExp3} are bounded by a multiple of $\Sjkzz$.

\begin{lma}
    \label{lma: CovNum}
    Suppose that assumptions A1, B1, B2, B4, and D1--D3 hold.  Then there exist constants $C_1$ and $C_2$ such that, for any $\epsilon > 0$, $(q,r) \in \{(0,0), (0,1), (1,0)\},$ and $j,k \in \{1,\ldots,p\},$
    \begin{equation*}
                \Pr\left(\normjk{\son \Wijk^{qr}} > \epsilon \right) \leq C_1\exp\left\{\frac{-C_2\epsilon^2}{\qjkn + \bgammaj\inv\bgammak\inv\vjkn\epsilon}\right\}.
    \end{equation*}
\end{lma}
\begin{proof}[\textbf{\upshape Proof:}]
    As the proof follows the same logic as that of Lemma~\ref{lma: MeanNum1}, the arguments will be sketched for the case $q = r = 0.$  For simplicity, write $\Wijk$ for $\Wijk^{00}.$  First, under assumptions D1 and D2, it can be deduced that, conditionally on $\mathbb{T}_{ijk},$ the random variables $\Uijklm$ are sub-Exponential random variables with parameters bounded by some universal constant $\rho^2 < \infty$ depending only on $\overline{\theta}$ and $\sigma$ in assumptions D1 and D2, respectively.  Hence, using moment bounds for sub-Exponential random variables, it can be concluded that, for any $j,k$, $(s,t) \in \Tj\times\Tk,$ and $\nu \geq 2,$
    \begin{equation*}
        E\left[\left|\Wijk(s,t)\right|^\nu \, | \,\mathbb{T}_{ijk}\right] \leq 2\nu!(2\rho^2\vijk)^\nu\left\{\sIijk \Kb{\bgammaj}(\Tijl - s)\Kb{\bgammak}(\Tikm - t)\right\}^\nu.
    \end{equation*}
    Then assumptions A1, B1, and B2 imply that
    $
    E\left[\left\{\sIijk \Kb{\bgammaj}(\Tijl - s)\Kb{\bgammak}(\Tikm - t)\right\}^2\right] \leq Ca_{ijk},
    $
    for some $C$, with
    $$
    a_{ijk} = \begin{cases}
        \Nij\Nik(\bgammaj\inv + \Nij - 1)(\bgammak\inv + \Nik - 1), & j \neq k, \\
        \Nij(\Nij - 1)\left\{\bgammaj\inv(\bgammaj\inv + \Nij - 2) + (\Nij - 2)(\bgammaj\inv + \Nij - 3)\right\}, & j = k,
    \end{cases}
    $$
    under an FR design, and $a_{ijk} \equiv a_i$ under an SR design, where $a_i$ is the common value of $a_{ijj}$ across $j$ in the above display.  Hence, applying Jensen's inequality, there are universal constants $B_1$ and $B_2$ such that
    \begin{equation*}
        \son E\left[\normjk{\Wijk}^\nu\right] \leq \frac{\nu!}{2}\left(B_2\son\vijk^2a_{ijk}\right)(B_1\bgammaj\inv\bgammak\inv\vjkn)^{\nu - 2},
    \end{equation*}
    where $\vjkn = \max_{i \in \{1,\ldots,n\}} \vijk|\Iijk|.$  Using the definitions of $\qjkn$ and applying Theorem 2.5 of \cite{bosq:00}, the inequality of the lemma is established for $q = r = 0$.  For the other values, apply the same arguments to the kernel $\tilde{K}(u) = uK(u).$
\end{proof}

\begin{lma}
    \label{lma: CovNumSup}
    Suppose that assumptions A1, B1, B2, and D1--D3 hold.  Then there exist constants $C_1,C_2,C_3 > 0$ such that, for any $\epsilon > 0,$ and $j,k \in \{1,\ldots,p\}$, 
    \begin{equation*}
        \Pr\left\{\left|\son \vijk \sIijk\left(|\Uijklm| - E\left[|\Uijklm|\right]\right)\right| > \epsilon\right\}  < C_1\exp\left\{\frac{-C_2\epsilon^2}{\son\vijk^2|\Iijk|^2 + \vjkn\epsilon}\right\}.
    \end{equation*}
\end{lma}

\begin{proof}
    The proof follows the same lines as the proof of Lemma~\ref{lma: MeanNum1Sup}, using the fact that each of $\Uijklm$ are, conditionally on $\mathbb{T}_{ijk}$, sub-Exponential random variables due to assumptions D1 and D2.  The details are omitted.
\end{proof}

%% file: CovTheoremsAndCorollariesJMVA.tex
\begin{proof}[\textbf{\upshape Proof of Theorem~\ref{thm: CovGen}:}]
The proof applies the same logic as the proof of Theorem~\ref{thm: meanGen}, but using Lemmas~\ref{lma: SjkqrConcIneq}--\ref{lma: CovNum} that are relevant for covariance estimation rather than the corresponding Lemmas~\ref{lma: SjrConcIneq}--\ref{lma: MeanNum1} that are for mean estimation.  Due to the similarities, step-by-step derivations of the bounds obtained below will not be provided.  As for mean estimation, the proof will be given for an FR design, noting that the same arguments can be followed under an SR design with weaker conditions on the bandwidth as in assumption D3 since $\Sjkqr$ do not depend on $j,k.$  Let $C_1,C_2,C_3$ be sufficiently large so that $
\limsup_{n \rightarrow \infty} \max_{j,k \in \{1,\ldots,p\}} \normsupjk{E[\Sjkzz]} \leq C_3
$ and Lemmas~\ref{lma: InfDenCov}--\ref{lma: CovNum} hold simultaneously.  Let $\omega_{jk}$ be as defined in \eqref{eq: omegajk}, and $S_{jk}^*(s,t)$, $\eta$ be as in Lemma~\ref{lma: InfDenCov}.  Recalling the expression of $\hgammajk(s,t) - \gammajk(s,t)$ in \eqref{eq: covDiff}, first apply Lemma~\ref{lma: InfDenCov} to conclude that, for $n \geq N$, 
$$
\Pr\left\{\min_{j,k \in \{1,\ldots,p\}} S_{jk}^*(s,t) < \eta\right\} \leq C_1p^{2 - C_2\eta^2\left\{\log(p)\max_{j,k \in \{1,\ldots,p\}} \bgammaj\inv\bgammak\inv \omega_{jk}^2\right\}}
$$
so that, by assumption D3, 
$
\max_{j,k \in \{1,\ldots,p\}} \left\{\inf_{(s,t) \in \Tj\times \Tk} S_{jk}^*(s,t)\right\}\inv = O_P(1).
$
Similarly, Lemma~\ref{lma: SupNumCov} implies that, with $Q_{jkr}$ as in \eqref{eq: Qjkr},
$
\max_{j,k \in \{1,\ldots,p\}} \normsupjk{Q_{jkr}} = O_P(1),
$
$r \in \{0,1,2\}$.  Hence, the rate is determined by that of the terms in \eqref{eq: CovNumTermExp1}--\eqref{eq: CovNumTermExp3}.  

Applying Lemma~\ref{lma: CovNum}, for any nonnegative integers $q$ and $r$ such that \mbox{$0 \leq q + r \leq 1,$} and any $R > 0,$
\begin{equation*}
    \Pr\left(\max_{j,k \in \{1,\ldots,p\}}\normjk{\Wijk^{qr}} > R\{\log(p)\}^{1/2}\max_{j,k \in \{1,\ldots,p\}}\left[\qjkn + \{\log(p)\}^{1/2}\bgammaj\inv\bgammak\inv\vjkn\right]\right) \leq C_1p^{2 - \frac{C_2R^2}{1 + R}}. 
\end{equation*}
Lemma~\ref{lma: SjkqrConcIneq} along with the bound $\normjk{\Sjkzz} \leq B\normsupjk{\Sjkzz}$ for a universal constant $B$ implies that, for any $R > 0,$
\begin{equation*}
        \Pr\left[\max_{j,k \in \{1,\ldots,p\}} \left(\normjk{\Sjkzz} - |\Tj|^{1/2}|\Tk|^{1/2}C_3\right) > R\{\log(p)\}^{1/2}\max_{j,k} \bgammaj\inv\bgammak\inv\omega_{jk}^2\right] \leq C_1p^{2 - C_2R^2}.
\end{equation*}
Since $\normsupjk{B_{jk\ell m}} \leq C\bgammaj\bgammak$ for a universal constant $C$, the proof is complete upon observing that
\begin{equation*}
\begin{split}
    &\max_{j,k \in \{1,\ldots,p\}}\Bigg\lVert\son \Wijk^{qr}(\cdot,*) + \son \vijk \sIijk \Kb{\bgammaj}(\Tijl - \cdot)\Kb{\bgammak}(\Tikm - *)\left(\frac{\Tijl - \cdot}{\bgammaj}\right)^q\left(\frac{\Tikm - *}{\bgammak}\right)^rB_{jklm}(\cdot,*)\Bigg\rVert_{j,k} \\
    &\hspace{1cm} = O_P\left(\max_{j \in \{1,\ldots,p\}} \bgammaj^2 + \max_{j,k \in \{1,\ldots,p\}} \left[\left\{\log(p)\qjkn\right\}^{1/2} + \log(p)\bgammaj\inv\bgammak\inv\vjkn \right]\right).
\end{split}    
\end{equation*}
\end{proof}

\begin{proof}[\textbf{\upshape Proofs of Corollary~\ref{cor: CovObsSubjL2}--\ref{cor: CovOptL2_sim}:}]
The proofs are similar to those of Corollaries~\ref{cor: meanObsSubjL2}--\ref{cor: meanOptL2_sim} and are omitted.
\end{proof}

\begin{proof}[\textbf{\upshape Proof of Theorem~\ref{thm: CovSupGen}:}]
The proof is similar to that of Theorem~\ref{thm: meanSupGen}, so only the details will be sketched.  First, as in the proof of Theorem~\ref{thm: CovGen}, it can be established that, with $S_{jk}^*$ as defined in Lemma~\ref{lma: InfDenCov} and $Q_{jkr}$ as defined in \eqref{eq: Qjkr} for $r \in \{0,1,2\}$, almost surely, $\max_{j,k \in \{1,\ldots,p\}} \left\{\inf_{(s,t) \in \Tj\times \Tk} S_{jk}^*(s,t)\right\}\inv = O(1)$ and $\max_{j,k \in \{1,\ldots,p\}} \normsupjk{Q_{jkr}} = O(1)$.
It then follows from \eqref{eq: covDiff}--\eqref{eq: CovNumTermExp3} that, almost surely,
\begin{equation}
    \label{eq: CovSupBnd}
    \max_{j,k \in \{1,\ldots,p\}} \normsupjk{\hgammajk - \gammajk} = O\left\{\max_{j \in \{1,\ldots,p\}}\bgammaj^2 +  \max_{j,k \in \{1,\ldots,p\}}\left(\normsupjk{\son\Wijk^{00}} + \normsupjk{\son\Wijk^{01}} + \normsupjk{\son\Wijk^{10}}\right)\right\}.
\end{equation}

Next, letting $\chi_{j}(\delta)$ represent an equally spaced grid of $\Tj$ with spacing no larger than $\delta > 0,$ set $\chi_{jk}(\delta) = \chi_j(\delta)\times \chi_k(\delta)$ such that $|\chi_{jk}(\delta)| \leq L\delta^{-2}$ for some universal constant $L$.  Then, for any nonnegative integers $q,r$ with $q + r \leq 1,$
\begin{equation*}
\normsupjk{\son \Wijk^{qr}} \leq \max_{(s,t) \in \chi_{jk}(\delta)} \left|\son \Wijk^{qr}(s,t)\right| + \sup_{|s - s'|,|t - t'| < \delta}\left|\son \left\{\Wijk^{qr}(s,t) - \Wijk^{qr}(s',t')\right\}\right|
\end{equation*}
Then, with $
a_{n}(\delta) = \max_{j,k \in \{1,\ldots,p\}}\left[\qjkn^{1/2} + \{\log(p\delta^{-2})\}^{1/2}\bgammaj\inv\bgammak\inv\vjkn\right],
$
using the arguments of Lemma~\ref{lma: CovNum}, it can be shown that,
for any $R > 0,$ there exist universal constants $C_1$ and $C_2$ such that
\[    
    \Pr\left[\max_{j,k \in \{1,\ldots,p\}}\max_{(s,t) \in \chi_{jk}(\delta)} \left|\son \Wijk^{qr}(s,t)\right| > R\left\{\log(p\delta^{-2})\right\}^{1/2}a_{n}(\delta) \right] \leq C_1(p\delta^{-2})^{1 - \frac{C_2R^2}{1 + R}},
\]
whence
\begin{equation}
    \label{eq: WijkqrGridMax}
\max_{j,k \in \{1,\ldots,p\}} \max_{(s,t) \in \chi_{jk}(n^{-\alpha})} \left|\son \Wijk^{qr}(s,t)\right| = O\left[\left\{\log(p\vee n)\qjkn\right\}^{1/2} + \log(p\vee n)\bgammaj\inv\bgammak\inv\vjkn\right]
\end{equation}
almost surely by Borel-Cantelli, where $\alpha$ satisfies assumption D4.  In addition, using assumption A2 and Lemma~\ref{lma: CovNumSup}, it follows that
\begin{equation}
    \label{eq: WijkqrGridDiff}
    \max_{j,k \in \{1,\ldots,p\}} \max_{|s - s'|,|t-t'| < \delta} \left|\son \left\{\Wijk^{qr}(s,t) - \Wijk^{qr}(s',t')\right\}\right| = O\left(\delta\max_{j \in \{1,\ldots,p\}} \bgammaj^{-3}\right)
\end{equation}
almost surely.  Applying $\delta = n^{-\alpha}$, $\alpha$ satisfying assumption D4, to \eqref{eq: CovSupBnd}--\eqref{eq: WijkqrGridDiff} proves the result.
\end{proof}

\begin{proof}[\textbf{\upshape Proofs of Corollaries~\ref{cor: CovObsSubjSup}--\ref{cor: CovOptSup_sim}:}] 
The proofs are similar to those of Corollaries~\ref{cor: CovObsSubjL2}--\ref{cor: CovOptL2_sim}, and are omitted.
\end{proof}